\documentclass[MSNbibl,number,citesort,dvips]{arxbj}
\usepackage{dcolumn}
\usepackage{multirow}
\usepackage{graphicx}

\aid{0}
\volume{18}
\issue{2}
\pubyear{2012}
\firstpage{606}
\lastpage{634}
\doi{10.3150/11-BEJ351}

\makeatletter
\newcolumntype{d}[1]{D{.}{.}{#1}}
\renewcommand{\emptyset}{\varnothing}
\newcommand{\cal}{\mathcal}
\newcommand{\argmin}{\arg\min}

\newtheorem{Example}{Example}
\newproclaim{Def}{Definition}
\newtheorem{Prop}{Proposition}
\newremark{Nota}{Remark}
\newtheorem{Lemm}{Lemma}
\newcommand{\fraca}[2]{{#1}/{#2}}
\makeatother

\begin{document}
\begin{frontmatter}

\title{Similarity of samples and trimming}
\runtitle{Similarity and trimming}

\begin{aug}
%%%% inicialai - be tarpu
\author[a]{\fnms{Pedro C.} \snm{\'Alvarez-Esteban}\corref{}\thanksref{a,e1}\ead[label=e1,mark]{pedroc@eio.uva.es}},
\author[a]{\fnms{Eustasio} \snm{del Barrio}\thanksref{a,e2}\ead[label=e2,mark]{tasio@eio.uva.es}},
\author[b]{\fnms{Juan~A.}~\snm{Cuesta-Albertos}\thanksref{b}\ead[label=e3]{juan.cuesta@unican.es}}
\and
\author[a]{\fnms{Carlos} \snm{Matr\'{a}n}\thanksref{a,e4}\ead[label=e4,mark]{matran@eio.uva.es}}
\runauthor{\'Alvarez-Esteban, del Barrio, Cuesta-Albertos and Matr\'{a}n}
\address[a]{Dept. de Estad\'\i stica e Investigaci\'on Operativa,
Universidad de Valladolid,
Prado de la Magdalena s.n.,
47005 VALLADOLID,
Spain.\\
\printead{e1,e2,e4}}
\address[b]{Dept. Matem\'{a}ticas, Estad\'\i stica y Computaci\'{o}n,
Universidad de Cantabria,
Avda. los Castros s.n.,
39005 SANTANDER,
Spain.
\printead{e3}}
\end{aug}

% HISTORY:
\received{\smonth{6} \syear{2010}}
\revised{\smonth{12} \syear{2010}}

% ABSTRACT
%
\begin{abstract}
We say that two probabilities are similar at level $\alpha$ if they are
contaminated versions (up to an $\alpha$ fraction) of the same common
probability.
We show how this model is related to minimal distances between sets of
\textit{trimmed probabilities}. Empirical versions turn out to present
an \textit{overfitting}
effect in the sense that trimming beyond the similarity level results
in trimmed samples that are closer than expected to each other. We show
how this can be
combined with a bootstrap approach to assess similarity from two data samples.
\end{abstract}

% KEYWORDS
%

\end{frontmatter}
%

%s1 ###
\section{Similarity vs. homogeneity}
Classical goodness of fit deals with the problem of assessing whether
the unknown random generator, $P$, of a data object, $X$, belongs to a
given class, $\mathcal{F}$. This includes two-sample problems in which
two different random objects are observed. We focus on checking whether
a certain feature of the corresponding random generators coincides. The
case in which $X_1$ is a collection of i.i.d. random variables
$X_1^1,\ldots,X_n^1$ with common distribution~$P_1$, $X_2$ is another
sequence of i.i.d. random variables $X^2_1,\ldots,X^2_m$ with law $P_2$
and the goal is to assess whether $\theta(P_1)=\theta(P_2)$ for some
function $\theta(\cdot)$ (including, for instance, \mbox{$\theta(P)=P$}) is a
\textit{homogeneity} problem, to which a large amount of literature has
been devoted. Our starting point is that it is often the case that the
researcher is not really
interested in checking whether $P\in\mathcal{F}$ or whether $P_1=P_2$.
Imagine the case of a~pharmaceutical company trying to introduce a new
(and cheaper) alternative to some
reference drug. The regulatory authorities will approve the new drug if
its performance with respect to a certain biological magnitude does not
differ from that of the standard drug. Both drugs could produce a
similar outcome on most patients. However, if there is a fraction of
them for whom the results are clearly different, then the new drug is
very likely to be rejected by a homogeneity test, while, in fact, it
has a similar performance for most individuals. As another example,
consider the comparison of two human populations that were initially
equal but have received immigration with different patterns. In these
situations the relevant assumption to check is not homogeneity, but
rather \textit{similarity} in the following sense.
\begin{Def}\label{asimilar}
Two probability measures $P_1$ and $P_2$ on the same sample space are
$\alpha$-si\-milar if there exist probability measures $P_0$, $P_1'$,
$P_2'$ such that
%
%e1 ###
\begin{equation}\label{asimilar2}
\cases{
P_1=(1-\varepsilon_1)P_0+\varepsilon_1 P_1',\cr
P_2=(1-\varepsilon_2)P_0+\varepsilon_2 P_2'}
\end{equation}
 with $0\leq\varepsilon_i \leq\alpha$, $i=1,2$.
\end{Def}

Definition \ref{asimilar} measures the overlap between $P_1$ and $P_2$,
in agreement with other possible
measures of similarity
(see the section ``Similarity between Populations'' in   \cite{Gower}).
Beware that smaller values of $\alpha$ in Definition \ref{asimilar}
correspond to more similar distributions (the case $\alpha=0$ being
equivalent to $P_1=P_2$).

A related situation, for one-sample problems, would be the case when we
observe some random object $X$ with law $P_1$. Ideally, $P_1$ should
equal $P_0$
(some gold standard), but the presence of noise means that, in fact,
%
%e2 ###
\begin{equation}\label{sim1lado}
P_1=(1-\varepsilon) P_0+ \varepsilon N,\qquad\varepsilon\leq\alpha
\end{equation}
for some unspecified $N$ if we assume that the noise level does not
exceed $\alpha$. We would say that
$P_1$ is \textit{similar} to $P_0$ at level $\alpha$ if (\ref
{sim1lado}) holds (observe that $P_1$ and $P_0$ do not
play a symmetric role in this definition).
In two-sample problems, we want to assess whether the two samples
can be assumed to be noisy realizations of some unkown gold standard,
as in Definition \ref{asimilar}.
Model (\ref{sim1lado}) corresponds to the `contamination neighborhoods'
introduced in Huber \cite{Huber64,Huber65} in a
robust testing setup. We discuss further connections to these and other
related references in  Section \ref{alternative}
below. Our goal in this work is to present a method for assessing
similarity of the unknown random generators
$P_1,P_2$ of two independent i.i.d. samples. Our procedure also yields
an estimate of the \textit{common core}
of the two distributions.

Our approach is based on trimming. Trimming procedures are of frequent
use in robust statistics as a way of downplaying the influence of
contaminating data in our inferences. The introduction of
data-dependent versions of trimming, often called impartial trimming,
allows us to overcome some limitations of earlier versions of trimming
that simply removed extreme observations at tails. Generally, impartial
trimming is based on some optimization criterion, keeping the fraction
of the sample (of a prescribed size) that yields the least possible
deviation with respect to a theoretical model. Today, impartial
trimming constitutes one of the main tools in the robust approach to a
variety of statistical settings
(see   \cite{Cuesta2,Garcia,Maronna2,Rousseeuw1}).
The first approach to model validation based on
impartial trimming is (to the best of our knowledge) the one in \'Alvarez-Esteban \textit{et al}. \cite{Pedro,Pedro3}. The problem considered there can be rephrased as
follows. Given two independent i.i.d. samples of univariate data with
unknown random generators $P_1,P_2$, we want to assess whether
$P_i=\mathcal{L}(\varphi_i (Z))$, $i=1,2$, for some random variable $Z$
defined on a probability space $(\Omega, \mathcal{F},\mathbb{P})$ and
non-decreasing functions, $\varphi_1$, $\varphi_2$, such that
\[
\mathbb{P}\bigl(\varphi_1(Z)\ne\varphi_2(Z)\bigr)\leq\alpha
\]
(see Section \ref{alternative} for further discussion). Despite the
interest of this approach, we believe that the similarity model given
by Definition \ref{asimilar} is often more natural and useful in
applications. Some technically related results and the connection with
the optimal transportation problem have been
reported in \'Alvarez-Esteban \textit{et al}. \cite{Pedro2}. A related approach based on density
estimation can be found in Mart\'inez-Camblor \textit{et al}. \cite{Camblor}.

As we will show in Section \ref{SecPrelim}, the similarity model of
Definition \ref{asimilar} can be expressed in terms of a minimal
distance between the sets of \textit{trimmings} of the probabilities
$P_i$, $i=1,2$. These are the sets of probabilities that one obtains
from a fixed one by removing or downplaying (to some degree) the weight
assigned by the original probability. When we look for the minimal
distance between trimmings of the empirical measures based on two
samples, we are
highlighting the part of the data that, hopefully, comes from the
common core $P_0$. From a descriptive point of view, this gives an
interesting tool for the comparison of data samples.

A distinctive feature of our proposal concerns the rates of
convergence. If $P_n$, $Q_n$ are the empirical distributions based on
two samples of univariate data (of equal size for simplicity), we will
trim up to an $\alpha$-fraction of data from both samples in order to
minimize some distance, $d(\cdot,\cdot)$; and if we write $P_{n,\alpha
}$, $Q_{n,\alpha}$ for the optimally trimmed empirical distributions,
we will have $d(P_{n,\alpha}, Q_{n,\alpha})\leq d(P_n, Q_n)$.
Trimming procedures generally give a balanced compromise between
efficiency and robustness, and increasing the level of trimming has a
moderate effect on the efficiency.
Thus, for univariate i.i.d. data coming from equal random generators,
we typically have $d(P_n, Q_n)=\mathrm{O}_P(n^{-1/2})$ and $d(P_{n,\alpha},
Q_{n,\alpha})=\mathrm{O}_P(n^{-1/2})$, but it is not true that $d(P_{n,\alpha},
Q_{n,\alpha})=\mathrm{o}_P(n^{-1/2})$
(see, e.g., Theorem A.1 in   \cite{Pedro}). However, for our procedure, over-trimming (i.e., trimming
beyond the similarity level) will produce an \textit{over-fitting effect},
namely, $d(P_{n,\alpha}, Q_{n,\alpha})=\mathrm{o}_P(n^{-1/2})$. That will be the
key for the statistical application of the procedure. Roughly speaking,
if two random samples are trimmed more than required to delete
contamination, then two samples far more similar than expected are
obtained and it is feasible to distinguish this pair of trimmed samples
from any other pair of non-trimmed, non-contaminated samples. We
formalize this idea in Section \ref{SecPrelim}. As in
\'Alvarez-Esteban \textit{et al}. \cite{Pedro},
our choice for the metric $d$ is the $L_2$ Wasserstein distance.

This over-fitting effect can be combined with a bootstrap procedure to
consistently decide if the underlying distributions of two i.i.d.
samples are similar in the sense of Definition \ref{asimilar} as we
show in Section \ref{bootstrapsection}. This statistical procedure
should also be useful in other frameworks of model validation. The
consistency of our procedure is independent of the kind of
contaminations. However, as expected, inliers are harder to detect than
outliers. In this proposal, we have to consider small resampling sizes
in the presence of inliers. This is discussed in Section \ref
{SecSimul}, where we present some simulations showing the performance
of our bootstrap procedure over finite samples. We also include the
analysis of a real data set.

For the sake of readability we have moved most of the proofs to an
\hyperref[appm]{Appendix}, together with some additional results on rates of convergence.

Throughout the paper $\mathcal{P}$ will be the set of Borel probability
measures on the real line,~$\mathbb{R}$, while ${\cal F}_p$ will denote
the set of distributions in ${\cal P}$ with finite $p$th absolute
moment. If $F$ is a distribution function, $F^{-1}$ will denote its
generalized inverse or quantile function.
Given $P,Q \in{\cal P}$, by $P \ll Q$ we will denote absolute
continuity of $P$ with respect to $Q$, and by $\frac{\mathrm{d}P}{\mathrm{d}Q}$ the
corresponding Radon--Nikodym derivative.
Unless otherwise stated, the random variables will be assumed to be
defined on the same probability space $(\Omega, \sigma, \nu)$. Weak
convergence of probabilities will be denoted by $ \to_{w} $ and ${\cal
L}(X)$ (resp., $EX$) will denote the law (resp., the
mean) of the variable $X$. The indicator function of a set $A$ will be
$I_{A}$ and $\ell$ will denote the Lebesgue measure.\vspace*{-1.5pt}

%s2 ###
\section{Trimming and over-fitting}\vspace*{-1.5pt} \label{SecPrelim}

%s2.1 ###
\subsection{Trimmings of a distribution}\vspace*{-1.5pt}

Trimming an $\alpha$-fraction of data in a sample of size $n$ can be
understood as replacing the empirical measure by a new one in which the
data are reweighted so that the trimmed points now have zero
probability while the remaining points will have weight $1/n(1-\alpha
)$. By analogy we can define the trimming of a distribution as follows.
\begin{Def}\label{DefiTrim}
Given $\alpha\in(0,1)$, we define the set of $\alpha$-trimmed versions
of $P$ by
%
%e3 ###
\begin{equation}\label{definition1}
\mathcal{R}_\alpha(P) := \biggl\{Q\in\mathcal{P}\dvt   Q\ll P, \frac
{\mathrm{d}Q}{\mathrm{d}P}\leq\frac{1}{1-\alpha},   P\mbox{-a.s.} \biggr\}.
\end{equation}
\end{Def}

This definition has been considered by several authors
(see  \cite{Gorda,Cascos2,Pedro}).
It allows the consideration of partial
removal of the points in the support of the probability. This
flexibility results in nice properties of the sets of trimmings, making
$\mathcal{R}_\alpha(P)$ a convex set, compact for the topology of weak
convergence (see Proposition 2.1 in   \cite{Pedro2}).

In this paper we use the quadratic Wasserstein distance, $\mathcal
{W}_2$, namely, the minimal quadratic transportation cost between
probabilities with finite second moment.~$\mathcal{W}_2$~met\-rizes weak
convergence plus convergence of second moments. We refer the reader to
Section 8 of Bickel and Freedman \cite{Bickel} for further details on $\mathcal W_{2}$. On
the real line $\mathcal W_{2}$ is simply the $L_2$ distance between
quantile functions, that is, $\mathcal{W}_2^2 ( P_1, P_2) =\int_{0}^1
(F_1^{-1}(t)-F_2^{-1}(t))^2\,\mathrm{d}t$ if $F_i^{-1}$ is the quantile function
of $P_i$. Trimmings are also well behaved with respect to $\mathcal
{W}_2$, as shown in \'Alvarez-Esteban \textit{et al}. \cite{Pedro2}. For instance, for $P \in\mathcal
{F}_2$, $\mathcal{R}_\alpha(P)$ is a compact subset of $\mathcal{F}_2$
for $\mathcal{W}_2$ (see Proposition 2.8 in   \cite{Pedro2}). A simple
consequence is that in
%
%e4 ###
\begin{equation}\label{minimaldist}
\mathcal{W}_2(\mathcal{R}_\alpha(P_1),\mathcal{R}_\alpha(P_2)):=\min
_{R_i\in\mathcal{R}_\alpha(P_i)}\mathcal{W}_2(R_1,R_2)
\end{equation}
the minimum is indeed attained. A remarkable result is that the
minimizer is unique under mild assumptions. This is Theorem 2.16 in
\'Alvarez-Esteban \textit{et al}. \cite{Pedro2}, which generalizes related results in
Caffarelli and McCann \cite{CaffarelliMcCann} and Figalli
\cite{Figalli}.\vadjust{\goodbreak}
\begin{Prop}\label{uniqueness}
If $P_1, P_2 \in\mathcal{F}_2$, $0<\alpha<1$ and $P_1$ or $P_2$ has a
density, then there exists a unique pair $(P_{1,\alpha},P_{2,\alpha
})\in\mathcal{R}_\alpha(P_1)\times\mathcal{R}_\alpha(P_2)$
such that
\[
\mathcal{W}_2(P_{1,\alpha},P_{2,\alpha})=\mathcal{W}_2(\mathcal
{R}_\alpha(P_1),\mathcal{R}_\alpha(P_2)),
\]
provided $\mathcal{W}_2(\mathcal{R}_\alpha(P_1),\mathcal{R}_\alpha(P_2))>0$.
\end{Prop}

The connection between trimmings and the similarity model of Definition
\ref{asimilar} is given by the next result. Here $d_{\mathrm{TV}}$ denotes the
distance in total variation, namely,
$d_{\mathrm{TV}}(P_1,P_2)=\sup_B |P_1(B)-P_2(B)|$, where $B$ ranges among all
Borel sets.
\begin{Prop}\label{equivalent} For $\alpha\in[0,1)$ the following are
equivalent:
\begin{enumerate}[(b)]
\item[(a)] $P_1$ and $P_2$ are $\alpha$-similar.
\item[(b)] $\mathcal{R}_\alpha(P_1)\cap\mathcal{R}_\alpha(P_2)\ne
\emptyset$.
\item[(c)] $d_{\mathrm{TV}}(P_1,P_2)\leq\alpha$.
\end{enumerate}
If $P_1,P_2\in\mathcal{F}_2$, then \textup{(a)}, \textup{(b)} or \textup{(c)} is equivalent to
\begin{enumerate}
\item[(d)] $\mathcal{W}_2 (\mathcal{R}_\alpha(P_1),\mathcal{R}_\alpha
(P_2))=0$.
\end{enumerate}
Finally, the common core distribution, $P_0$, in Definition \ref
{asimilar} is unique if and only if \mbox{$d_{\mathrm{TV}}(P_1,P_2)=\alpha$}.
In this case, $P_0$ is given by the density $f_0=(f_1\wedge
f_2)/(1-\alpha)$ with respect to $\mu$ if
$\mu$ is a common $\sigma$-finite dominating measure for $P_1$ and
$P_2$ and $f_1$ and $f_2$ are the corresponding densities and we have the
canonical decomposition $P_i=(1-\alpha)P_0+\alpha P_i'$, $i=1,2$,
$P_i'$ having density $\frac1 {\alpha} (f_i-f_1\wedge f_2)$ with
respect to $\mu$.
\end{Prop}

\begin{pf} If (a) holds, then $P_0(A)\leq\frac{1}
{1-\alpha}P_i(A)$ for all Borel $A$. In particular, $P_0\ll P_i$ and,
if $A_i=\{\frac{\mathrm{d}P_0}{\mathrm{d}P_i} > (1-\alpha)^{-1}\}$, obviously
$P_0(A_i)=0$ and $P_0\in\mathcal{R}_\alpha(P_1)\cap\mathcal
{R}_\alpha(P_2)$, showing (b). Assume now (b) and take $P_0\in\mathcal
{R}_\alpha(P_1)\cap\mathcal{R}_\alpha(P_2)$. Then $(1-\alpha
)P_0(A)\leq P_i(A)$ for all $A$. If $\alpha=0,$ then (c) holds
trivially. Otherwise define $P_i'(A)=(P_i(A)-(1-\alpha)P_0(A))/\alpha$.
Then~$P_i'$ is a probability and $d_{\mathrm{TV}}(P_1,P_2)=\alpha
d_{\mathrm{TV}}(P_1',P_2')\leq\alpha$, that is, (c) holds. Finally, we assume
that (c) holds and take
$\mu$ to be a common $\sigma$-finite dominating measure for~$P_1$ and~$P_2$
and write $f_1$ and $f_2$ for the corresponding densities.
Then (see Lemma 2.20 in~\cite{Massart}) $d_{\mathrm{TV}}(P_1,P_2)=1-\int(f_1
\wedge f_2)\,\mathrm{d}\mu$ (where $a\wedge b$ means $\min(a,b)$).
Write $\varepsilon=d_{\mathrm{TV}}(P_1,P_2)$ and assume $\varepsilon>0$ (the
case $\varepsilon=0$ is trivial).
We set $f_i'=(f_i-f_1 \wedge f_2)/\varepsilon$, $i=1,2$, and
$f_0=(f_1\wedge f_2)/(1-\varepsilon)$. $f_0,f_1',f_2'$
are densities with respect to $\mu$. We write $P_0, P_1', P_2'$ for the
associated probabilities.
Then (\ref{asimilar}) holds with $\varepsilon_1=\varepsilon
_2=\varepsilon\leq\alpha$.
Equivalence of~(b) and~(d) follows from compactness of the sets of
trimmings. The last claim follows easily from the arguments above.
\end{pf}

\begin{Nota}\label{desccanonica} It follows from Proposition \ref
{equivalent} that $\mathcal{W}_2 (\mathcal{R}_\alpha(P_1),\mathcal
{R}_\alpha(P_2)) >0$ if and only if $d_{\mathrm{TV}}(P_1,P_2)>\alpha$, that
is, $d_{\mathrm{TV}}(P_1,P_2)$ is the minimal level of trimming required to make
$P_1$ and $P_2$ equal. Also, if $d_{\mathrm{TV}}(P_1,P_2)=\alpha$, then
the probability $P_0$ with density $f_0=(f_1\wedge f_2)/(1-\alpha)$
with respect to $\mu$ (as in the proof above) is the unique element
in\vadjust{\goodbreak}
$\mathcal{R}_\alpha(P_1)\cap\mathcal{R}_\alpha(P_2)$. This means
that, as in Proposition \ref{uniqueness}, there is also a unique pair,
namely, $(P_0,P_0)\in\mathcal{R}_\alpha(P_1)\times\mathcal{R}_\alpha
(P_2)$ such that
\[
\mathcal{W}_2(P_0,P_0)=\mathcal{W}_2(\mathcal{R}_\alpha(P_1),\mathcal
{R}_\alpha(P_2))=0.
\]
This extends the result in Proposition \ref{uniqueness} to the case
$d_{\mathrm{TV}}(P_1,P_2)\geq\alpha$.
\end{Nota}

Proposition \ref{equivalent} shows that the similarity model (\ref
{asimilar}) can be expressed in terms of different metrics.
In fact, (d) would remain true if $\mathcal{W}_2$ were replaced by any
other metric for which the sets of trimmings are compact.
With applications in mind, $\mathcal{W}_2$ turns out to be a more
convenient choice. In order to assess (\ref{asimilar}) from two samples
of i.i.d. data with empirical distributions $P_{1,n}$ and $P_{2,m}$,
say, we will have $d_{\mathrm{TV}}(P_{1,n},P_{2,m})=1$ almost surely (provided~$P_1$
and $P_2$ have densities)
and we cannot use (at least in a na\"ive fashion) formulation~(c). On
the other hand, $\mathcal{W}_2$ is well behaved in this respect and
empirical versions of both the minimal distances and the minimizers are
consistent estimators of their theoretical counterparts. This is the
content of the following result (Theorem 2.17 in
  \cite{Pedro2}). We quote it here for completeness.

\begin{thm}[(Consistency)]\label{consistencia}
Let $ \lbrace X_{n} \rbrace_n$, $ \lbrace Y_{n}
\rbrace_n$ be two sequences of i.i.d. random variables with $\mathcal
L(X_n)=P$, $\mathcal L(Y_n)=Q$, $P, Q \in\mathcal{F}_2$, and write
$P_{n}$, $Q_{m}$ for the empirical distributions based on the samples
$X_1,\ldots,X_n$ and $Y_1,\ldots, Y_m$, respectively. Then,
if $\min(m,n)\to\infty$,
\[
\mathcal{W}_2(\mathcal{R}_\alpha(P_{n}),\mathcal{R}_\alpha(Q_{m}))\to
\mathcal{W}_2(\mathcal{R}_\alpha(P),\mathcal{R}_\alpha(Q))\qquad\mbox{a.s.}
\]
Further, if $P$ or $Q \ll\ell$ and $d_{\mathrm{TV}}(P,Q)\geq\alpha,$ then
\[
\mathcal W_2(P_{n,\alpha} , P_{\alpha}) \to0  \quad \mbox{and} \quad  \mathcal
W_2(Q_{m,\alpha} , Q_{\alpha}) \to0\qquad\mbox{a.s.},
\]
where $(P_{\alpha},Q_{\alpha})=\argmin_{R_1\in\mathcal
{R}_\alpha(P),R_2\in\mathcal{R}_\alpha(Q)} \mathcal{W}_2(R_1,R_2)$ and
$(P_{n,\alpha},Q_{m,\alpha})$ are defined similarly from
$P_{n}$, $Q_{m}$.
\end{thm}

%s2.2 ###
\subsection{Related concepts and works}\label{alternative}

The similarity model (\ref{asimilar2}) is obviously related to the
so-called `contamination neighborhoods' of a probability $P_0$, defined as
%
%e5 ###
\begin{eqnarray}\label{contamination}
\mathcal{V}_\varepsilon(P_0)&:=&\{(1-\varepsilon)P_0+\varepsilon P'\dvt
P' \in\mathcal{P}\} \nonumber
\\[-8pt]
\\[-8pt]
&\hspace*{2.6pt}=& \{Q\in\mathcal{P}\dvt Q(A)\leq(1-\varepsilon
)P_0(A)+\varepsilon\mbox{ for every Borel set } A \},
\nonumber
\end{eqnarray}
which have been widely used in the theory of robust statistics after
the pioneering works by Huber \cite{Huber64,Huber65}. In particular,
Huber \cite{Huber65} introduced these
neighborhoods in robust testing, providing a robust version of the
Neyman--Pearson lemma for simple hypothesis versus simple alternative.
This theory was completed for more general sets of hypotheses and
alternatives, additionally considering more flexible neighborhoods in
Huber and Strassen \cite{HuberStrassen}, Rieder \cite{Rieder}
and Buja \cite{Buja}. In fact, Rieder's neighborhoods of a probability~$P_0$, defined as
%
%e6 ###
\begin{equation}\label{Riederneigh}
\mathcal{V}^R_{\varepsilon,\delta} (P_0):=\{Q\in\mathcal{P}\dvt Q(A)\leq
(1-\varepsilon)P_0(A)+\varepsilon+\delta\mbox{ for every Borel set }
A \},
\end{equation}
comprise contamination as well as total variation norm neighborhoods
(taking $\delta=0$ or $\varepsilon=0$, resp.).

It can be easily shown (see also Proposition 2.1 in   \cite{Pedro2})
that $P\in\mathcal{V}_\varepsilon(P_0)$ is equivalent to $P_0\in
\mathcal{R}_\varepsilon(P)$.
Thus, our statement \textit{$P_1$ and $P_2$
are $\alpha$-similar} can also be expressed, in terms of contamination
neighborhoods, as \textit{there exists a probability $P_0$ such
that $P_1,P_2 \in\mathcal{V}_\alpha(P_0)$}. However, there are
different possibilities for such $P_0$, and the model
considered in this paper, given through any one of the equivalent
statements in Proposition \ref{equivalent}, cannot be expressed
in terms of a neighborhood, like (\ref{contamination}) or (\ref
{Riederneigh}) of a fixed probability.

\begin{figure}

\includegraphics{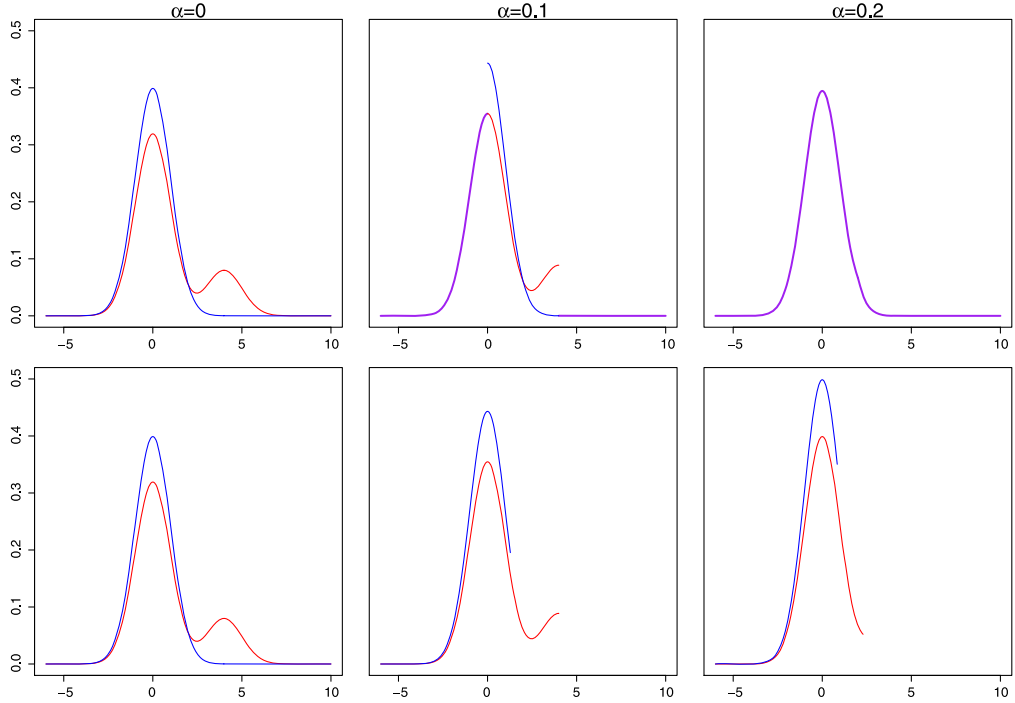}

\caption{Densities of optimally trimmed $P$ and $Q$ with independent
trimming (first row) and common trimming (second row).}\vspace*{6pt}
\label{FigNueva}
\end{figure}

Further related work includes \'Alvarez-Esteban \textit{et al}. \cite{Pedro}, where it is shown, for a
probability, $P$, on the real line, that $\mathcal{R}_\alpha(P)$ can be
expressed in terms of the trimmings of the uniform law on $(0,1)$,
$U(0,1)$. This set can be identified with the set $\mathcal{C}_{\alpha
}$ of absolutely continuous functions $h\dvtx[0,1] \rightarrow[0,1]$ such
that $h(0)=0$, $ h(1)=1$, with derivative $h'$ such that $0 \le h' \le
\frac{1}{1-\alpha}$. For function $h$, it is useful to write $P_h$ for
the probability measure with distribution function $h(P(-\infty,t])$. Then
%
%e7 ###
\begin{equation}\label{caparameter}
\mathcal{R}_\alpha(P)=\{P_h\dvt  h\in\mathcal{C}_\alpha\}.
\end{equation}
Hence, we can measure the deviation between the sets of trimmings of
$P$ and $Q$ through
\[
\mathcal{T}_\alpha(P,Q):=\min_{h\in\mathcal{C}_\alpha}\mathcal{W}_2(P_h,Q_h).
\]
We call $\mathcal{T}_\alpha(P,Q)$ the \textit{common trimming} distance
between $P$ and $Q$. If $P$ and $Q$ have quantile functions $F^{-1}$
and $G^{-1}$, then a simple change of variable shows
\begin{eqnarray*}
\mathcal{W}_2(P_h,Q_h)&=&\int_0^1
\bigl(F^{-1}(h^{-1}(x))-G^{-1}(h^{-1}(x))\bigr)^2\,\mathrm{d}x\\
&=&\int_0^1 \bigl(F^{-1}(y)-G^{-1}(y)\bigr)^2 h'(y)\,\mathrm{d}y.
\end{eqnarray*}
Thus, $\mathcal{T}_\alpha(P,Q)=0$ if and only if $\ell(\{y\in(0,1)\dvt
F^{-1}(y)\ne G^{-1}(y) \})\leq\alpha$. It follows easily from this
that $\mathcal{T}_\alpha(P,Q)=0$ if and only if there is a random
variable $Z$ defined on a~probability space $(\Omega, \mathcal
{F},\mathbb{P})$ and non-decreasing, left-continuous functions, $\varphi
_1$, $\varphi_2$, with $\mathcal{L}(\varphi_1(Z))=P$, $\mathcal
{L}(\varphi_2(Z))=Q$ such that
%
%e8 ###
\begin{equation}\label{corruptsignal}
\mathbb{P}\bigl(\varphi_1(Z)\ne\varphi_2(Z)\bigr)\leq\alpha.
\end{equation}
In contrast, since $d_{\mathrm{TV}}(P,Q)=\min\{\mathbb{P}(X\ne Y)\dvt  \mathcal
{L}(X)=P, \mathcal{L}(Y)=Q \}$ (see Lemma 2.20 in~\cite{Massart}), we
see that $\mathcal{W}_2(\mathcal{R}_\alpha(P),\mathcal{R}_\alpha(Q))=0$
if and only if $\mathcal{L}(\varphi_1(Z))=P$, $\mathcal{L}(\varphi
_2(Z))=Q$ for some random
variable $Z$ and measurable (not necessarily monotonic) $\varphi_i$
such that~(\ref{corruptsignal}) holds. In summary, \textit{two random
objects are $\alpha$-similar
if and only if they are different transforms of a common random signal
and the transforms differ from each other with probability at most
$\alpha$};
they are equivalent in terms of common trimming if and only if they are
different \textit{monotonic} transforms of a common random signal and the
transforms differ from each other with probability at most $\alpha$. In
the somewhat artificial event that we believe that our two samples come
from a monotonic, possibly different, transform of some original
signal, then the common trimming similarity model is reasonable.
Otherwise, the similarity model (\ref{asimilar}) is the natural choice.
For a less technical illustration of this idea we show in Figure \ref
{FigNueva} the different effect of independent and common trimming.
We have taken $P=N(0,1)$, $Q=0.8N(0,1)+0.2N(4,1)$ and three values of
the trimming level,~$\alpha$.
In the first row we show the densities of $P_\alpha$ (blue line) and
$Q_\alpha$ (red line), with $(P_\alpha, Q_\alpha)=\argmin_{R_1\in
\mathcal{R}_\alpha(P),
R_2\in\mathcal{R}_\alpha(Q)} \mathcal{W}_2(R_1,R_2)$. In this case,
trimming $\alpha=0.2$ results in $P_\alpha=Q_\alpha$, that is, trimming
removes contamination.
The second row shows the densities of $P_{h_\alpha}$ (blue line) and
$Q_{h_\alpha}$ (red line), with $h_\alpha=\argmin _{h\in\mathcal
{C}_\alpha}\mathcal{W}_2(P_h,Q_h)$.
Clearly, $P_{h_\alpha}$ and $Q_{h_\alpha}$ are different and this
remains true no matter how close to 1 we choose $\alpha$. If
trimming is used with the goal of
removing contamination and assessing that the core of the two
distributions are equal, then it is clear that the common trimming
approach fails to do so.

%f1 ###

In \'Alvarez-Esteban \textit{et al}. \cite{Pedro3} we have considered, under this common trimming setup,
the problem of testing
whether a random sample can be considered `mostly normal', that is, if
the generator of the sample is similar to a normal distribution with
unknown parameters.

Finally let us mention the application in \'Alvarez-Esteban \textit{et al}. \cite{Pedro2} of some
asymptotic results for a related two-sample problem:
Given
$X_1,\ldots,X_n$ i.i.d. $P$ and $Y_1,\ldots,Y_m$ i.i.d.~$Q,$
we consider testing the related null hypotheses
\begin{eqnarray*}
H_1\dvt  \mathcal{W}_2(\mathcal{R}_\alpha(P), \mathcal{R}_\alpha
(Q))&\leq& \Delta_0 \quad\mbox{vs.}\quad\mathcal{W}_2(\mathcal
{R}_\alpha(P), \mathcal{R}_\alpha(Q))> \Delta_0,
\\
H_2\dvt  \mathcal{W}_2(\mathcal{R}_\alpha(P), \mathcal{R}_\alpha
(Q))&\geq& \Delta_0 \quad\mbox{vs.}\quad\mathcal{W}_2(\mathcal
{R}_\alpha(P), \mathcal{R}_\alpha(Q))< \Delta_0
\end{eqnarray*}
for a given threshold $\Delta_0>0$ to be chosen by the practitioner.
Observe that rejecting the null hypothesis $H_2$ allows us to conclude
that, with high confidence, the unknown
random generators $P$ and $Q$ are not far from similarity.

%s2.3 ###
\subsection{The over-fitting effect of trimming}\label{overfitting}

In this subsection we keep the notation of Theorem \ref{consistencia}
and assume that we deal with two independent samples, $X_1,\ldots,X_n$
i.i.d. $P$ and
$Y_1,\ldots,Y_m$ i.i.d. $Q$. We write $P_n$, $Q_m$ for the empirical
measures and $P_{n,\alpha}$, $Q_{m,\alpha}$ are minimizers of the
$\mathcal{W}_2$ distance between trimmings of the empirical
distributions $P_n, Q_m$.

It follows from Theorem \ref{consistencia} that $\mathcal
{W}_2(P_{n,\alpha}, Q_{m,\alpha})\to0$ a.s. when the similarity model
(\ref{asimilar}) holds true and we may wonder about the rate of
convergence in this limit. Note that under homogeneity, that is, if
$P=Q$ and taking $n=m$ for simplicity, we have under integrability assumptions
%
%e9 ###
\begin{equation}\label{weaklimit}
\sqrt{n} \mathcal{W}_2(P_{n}, Q_{n})\to_w \biggl (2 \int_0^1 \frac
{B^2(t)}{f^2(F^{-1}(t))}\,\mathrm{d}t  \biggr)^{1/2},
\end{equation}
where $B$ is a Brownian bridge and $f$ and $F^{-1}$ are the density and
quantile functions of $P$
(this follows easily, for instance, from Theorem 4.6 in  \cite{Barrio}).
Thus, random samples from homogeneous generators have empirical
distributions at $\mathcal{W}_2$-distance of exact order~$n^{-1/2}$,
while, for non-homogeneous random generators
$\mathcal{W}_2(P_{n}, Q_{n})\to\mathcal{W}_2(P, Q)$, a positive
constant. Likewise, in the common trimming model of Section \ref
{alternative}, if $h_{n,\alpha}$ is such that
$\mathcal{T}_\alpha(P_n,Q_n)=\mathcal{W}_2((P_n)_{h_{n,\alpha}},
(Q_n)_{h_{n,\alpha}})$ and we write $\tilde{P}_{n,\alpha
}=(P_n)_{h_{n,\alpha}}$, $\tilde{Q}_{n,\alpha}=(Q_n)_{h_{n,\alpha}}$
(the optimal trimmings of the empirical measures), then, under $\mathcal
{T}_\alpha(P,Q)=0,$ we have that
$\sqrt{n}\mathcal{W}_2(\tilde{P}_{n,\alpha}, \tilde{Q}_{n,\alpha})$
converges in law to a non-null limit (Theorem A.1 in~\cite{Pedro}),
whereas if $\mathcal{T}_\alpha(P,Q)>0,$ then
$\mathcal{W}_2(\tilde{P}_{n,\alpha}, \tilde{Q}_{n,\alpha})$ converges
a.s. to a positive constant.

In the similarity model (\ref{asimilar}) the gap between the null and
the alternative is of higher order. If $P$ and $Q$ are not similar at
level $\alpha$, then $\mathcal{W}_2(P_{n,\alpha}, Q_{m,\alpha})\to
\mathcal{W}_2(P_\alpha, Q_\alpha)>0$ (Theorem \ref{consistencia}). On
the other hand, if $d_{\mathrm{TV}}(P,Q)<\alpha$, then our next result shows
that $\sqrt{n} \mathcal{W}_2(P_{n,\alpha}, Q_{n,\alpha})\to0$ in probability.
To avoid integrability issues, we assume $P$ and $Q$ to have bounded support;
this is enough for applications, since a monotonic transformation of the
data could achieve boundedness while preserving the distance in total
variation. Furthermore, it ensures that the conditions
$d_{\mathrm{TV}}(P,Q)\leq\alpha$ and $\mathcal{W}_2 (\mathcal{R}_\alpha
(P),\mathcal{R}_\alpha(Q))=0$ are equivalent.
\begin{thm}\label{rate}
Assume $P,Q \in\mathcal F_2$ are supported in a common bounded
interval and have densities bounded away from zero and with bounded
derivatives. Assume further that $n/(n+m)\to\lambda\in(0,1)$. If $\alpha
_n\in(0,1)$ satisfies $\alpha_n \geq d_{\mathrm{TV}}(P,Q) + \frac{r_n}{\sqrt
{n}}$ for some $r_n\to\infty$, then
%
%e10 ###
\begin{equation}\label{rate1}
\sqrt{n} \mathcal{W}_2(P_{n,\alpha_n},Q_{m,\alpha_n}) \to0   \qquad \mbox{in
probability.}
\end{equation}
\end{thm}

\begin{figure}

\includegraphics{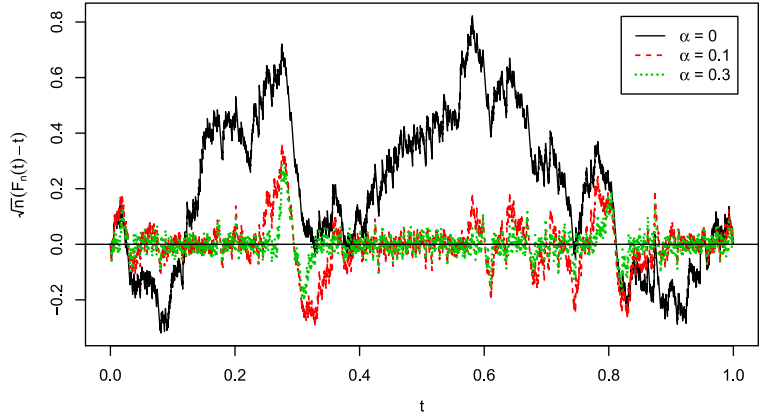}

\caption{Trajectories of the uniform empirical process (solid line) and two
variants based on trimming. The trimming levels are $\alpha=0.1$ and
$\alpha=0.3$ (dashed and dotted lines).}
\label{grafico7}
\end{figure}

We give a proof of Theorem \ref{rate} in the \hyperref[appm]{Appendix}. A similar
over-fitting effect is observed if a sample is over-trimmed to
optimally fit a given model: If $X_1,\ldots,X_n$ are i.i.d. $P$,
$P_{n,\alpha}=\argmin_{R\in\mathcal{R}_\alpha(P_n)}\mathcal{W}_2
(R,Q)$ and $\mathcal{W}_2 (\mathcal{R}_{\alpha_0}(P),Q)=0$ for some
$\alpha_0<\alpha$, then (see Theorem \ref{unamuestra} in the \hyperref[appm]{Appendix})
\[
\sqrt{n} \mathcal{W}_2(P_{n,\alpha},Q) \to0  \qquad  \mbox{in probability.}
\]
Empirical evidence of this over-fitting effect is shown in Figure \ref
{grafico7}. A random sample of size $n=1000$ from a $U(0,1)$
distribution was taken. This sample was trimmed using the proportions
$\alpha=0, 0.1, 0.3$ in order to obtain a sample as close to the
$U(0,1)$ as possible. We denote by $F_n^\alpha$ the distribution
function of $P_{n,\alpha}$ and in Figure \ref{grafico7}, we represent
the empirical processes $D_n^\alpha(t)=n^{1/2}(F_n^\alpha(t) -t)$, $t
\in[0,1]$ for $\alpha=0, 0.1, 0.3$.

%f2 ###

Since the true random generator and the target are the same, no
trimming is required in this case to remove contamination and, for
$\alpha>0,$ we are over-trimming.
Observe that~$D_n^{0.1}$ and $D_n^{0.3}$ do not differ too much from
each other, while they are quite far from the untrimmed version.

%s3 ###
\section{A bootstrap assessment of similarity}\label{bootstrapsection}

We show in this section how we can use the over-fitting effect of
trimming for the assessment of the similarity model (\ref{asimilar2}).
Again, we will assume that we observe two independent random samples
$X_1,\ldots,X_n$ i.i.d. $P$, $Y_1,\ldots,Y_m$ i.i.d. $Q$. We
would like to test the null hypothesis $H_0\dvt  d_{\mathrm{TV}}(P,Q)\leq\alpha$.
Theorem \ref{rate} says that trimming beyond the
similarity level kills randomness and results in (trimmed) samples that are
more similar to each other than random samples coming from the same generator.
We will use a bootstrap approach to generate suitable random samples
from a common
generator and compare the optimally trimmed distance to the distance
computed on the bootstrap replicates.

We write
$P_{n}$, $Q_{m}$ for the empirical distributions and, given $\alpha
_n\in(0,1)$,
\[
(P_{n,\alpha_n},Q_{m,\alpha_n})= \operatorname{\arg\min}\limits_{R_1 \in\mathcal{R}_{{\alpha
}_n}(P_n), R_2 \in\mathcal{R}_{\alpha_n}(Q_m)}\mathcal
{W}_2(R_1, R_2),
\]
so that $\mathcal{W}_2(P_{n,\alpha_n},Q_{m,\alpha_n})=\mathcal
{W}_2(\mathcal{R}_{\alpha_n}(P_n),\mathcal{R}_{\alpha_n}(Q_m))$.

We consider now the pooled probability
\[
R_{n,m} =\frac n {n+m}P_{n,\alpha_n}+ \frac m {n+m} Q_{m,\alpha_n}.
\]
$R_{n,m}$ is a random probability measure concentrated on $\{
Z_1,\ldots,Z_{n+m}\},$ where $Z_{j}=X_{j}$ for $j=1,\ldots,n$, and
$Z_{j}=Y_{j-n}$ for $j=n+1,\ldots,n+m$.

Conditionally, given the data, we draw new random variables,
$X_1^*,\ldots,X_{n'}^*, Y_1^*,\ldots,Y_{m'}^*$ i.i.d. $R_{n,m}$,
with $m'=[n' m/n]$ and $n'$ to be chosen later. We will use the
notation $\mathbb{P}^*$ for the bootstrap probability,
that is, the conditional probability given the original data $\{X_n\}
_n$, $\{Y_m \}_m$. Finally, by $P_{n'}^{*}$ and $Q_{m'}^{*}$
we will denote the empirical measures based on $X_1^*,\ldots,X_{n'}^*$ and
$Y_1^*,\ldots,Y_{m'}^*$, respectively. Now, we define
%
%e11 ###
\begin{equation}\label{bootstrappvalue}
p_{n,m}^*:=\mathbb{P}^*  \Biggl\{ \sqrt{\frac{n'm'}{n'+m'}} \mathcal
{W}_2(P_{n'}^{*},Q_{m'}^{*}) > \sqrt{\frac{nm}{n+m}} \mathcal
{W}_2(P_{n,\alpha_n},Q_{m,\alpha_n}) \Biggr\}.
\end{equation}
$p_{n,m}^*$ is the bootstrap $p$-value for the similarity model (\ref
{asimilar2}), with rejection for small values of it.
In practice $p_{n,m}^*$ can be approximated by Monte Carlo simulation.
We note that
if $n\alpha_n$ and $m\alpha_n$ are integer, typically the trimming
process will not produce partially trimmed points and $P_{n,\alpha_n}$
and $Q_{m,\alpha_n}$ will be the empirical measures on the sets of
non-trimmed data. If we take $\alpha_n\to\alpha$, then if the
similarity model fails, $\mathcal{W}_2(P_{n,\alpha_n},Q_{m,\alpha_n})$
will be large while $\mathcal W_2(P_{n'}^{*},Q_{m'}^{*})$ will vanish.
On the other hand, for similar distributions $\mathcal{W}_2(P_{n,\alpha
_n},Q_{m,\alpha_n})$ will vanish at a faster
rate than $\mathcal W_2(P_{n'}^{*},Q_{m'}^{*})$ and rejection for small
bootstrap $p$-values will result in a consistent rule. We make this
precise in our next result.

\begin{thm}\label{principal}
With the above notation, assume that $P,Q$ have densities satisfying
the assumptions of Theorem \ref{rate}. Assume further that
$n/(n+m)\to\lambda\in(0,1)$ and take $\alpha_n=\alpha+K/{\sqrt{n\wedge
m}}$ with $K>0$. Then, if $n'\to\infty$ and $n'=\mathrm{O}(n)$,
\begin{enumerate}[(ii)]
\item[(i)] if $d_{\mathrm{TV}}(P,Q)<\alpha,$ then $p_{n,m}^*\to1$ in
probability,
\item[(ii)] if $d_{\mathrm{TV}}(P,Q)>\alpha,$ then $p_{n,m}^*\to0$ in
probability.\vadjust{\goodbreak}
\end{enumerate}
\end{thm}

A proof of Theorem \ref{principal} is given in the \hyperref[appm]{Appendix}. It roughly
says that a test of the similarity model
(\ref{asimilar2}) that
rejects $\alpha$-similarity for values of $p_{n,m}^*$ above a fixed
threshold $L\in(0,1)$ is a consistent rule.
In order to make a sensible choice of the threshold, $L$, as well as of
the constant, $K$, in Theorem \ref{principal},
we still need to control the probability of rejection
at the boundary of the null hypothesis; that is, in the case
$d_{\mathrm{TV}}(P,Q)=\alpha$.
In this case we write again $P_0$ for the common part of $P$ and $Q$ in
the canonical decomposition
in Remark \ref{desccanonica}.
If $\tilde{P}_n\in\mathcal{R}_{\alpha_n}(P)$ and
$\tilde{Q}_n\in\mathcal{R}_{\alpha_n}(Q)$, with $\alpha_n$ as in
Theorem \ref{principal},
are such that $\mathcal{W}_2(\tilde{P}_n,\tilde{Q}_n)\to0,$ then, by
uniqueness, we have
$\mathcal{W}_2(\tilde{P}_n,P_0)\to0$. We introduce the
following assumption about rates in this convergence: If $\tilde{P}_n\in
\mathcal{R}_{\alpha_n}(P)$,
$\tilde{Q}_n\in\mathcal{R}_{\alpha_n}(Q)$ (and $\alpha
_n=d_{\mathrm{TV}}(P,Q)+\frac{K}{\sqrt{n}}$), then,
for some $\rho\in(0,1],$
%
%e12 ###
\begin{equation}\label{extrahipotesis}
\mathcal{W}_2(\tilde{P}_n,\tilde{Q}_n)=\mathrm{O}(n^{-1/2}) \quad \Rightarrow \quad \mathcal
{W}_2(\tilde{P}_n,P_0)=\mathrm{O}(n^{-\rho/2}).
\end{equation}
Under this assumption we can control the type I error probability using
our next result.

\begin{thm}\label{controldelnivel}
Under the assumptions and notation of Theorem \ref{principal}, if $P$
and $Q$ are such that
$d_{\mathrm{TV}}(P,Q)=\alpha$ and satisfy (\ref{extrahipotesis}), taking $n'\to
\infty$, $n'=\mathrm{o}(n^\rho)$ and
\[
\alpha_n=\alpha+\frac{\sqrt{\alpha(1-\alpha)}}{\sqrt{n\wedge m}}\Phi^{-1}
\bigl(\sqrt{1-\gamma}\bigr)
\]
with $\gamma\in(0,1)$, then
$\limsup_n \mathbb{P}(p_{n,m}^*\leq\beta)\leq\beta+\gamma$.
\end{thm}

The main consequence is that we can test the similarity model (\ref
{asimilar2}) at a given level $\beta+\gamma\in(0,1)$.
To be precise, if we replace our ideal $H_0\dvt  d_{\mathrm{TV}}(P,Q)\leq\alpha$
by $\tilde{H}_0$ consisting of pairs $(P,Q)$
satisfying the assumptions in Theorem \ref{principal} and
$d_{\mathrm{TV}}(P,Q)<\alpha$ or $d_{\mathrm{TV}}(P,Q)=\alpha$
plus Condition (\ref{extrahipotesis}), then, if
we reject for $p_{n,m}^*\leq\beta$, Theorems \ref{principal} and \ref
{controldelnivel} ensure
\[
\sup_{(P,Q)\in\tilde{H}_0} \limsup_n \mathbb{P}_{(P,Q)}(p_{n,m}^*\leq
\beta)\leq\beta+\gamma,
\]
where $\mathbb{P}_{(P,Q)}$ denotes probability assuming the laws of the
$X$'s and the $Y$'s are $P$ and $Q$,
respectively. It is in this sense that we can say that the procedure is
conservative, having an asymptotic level of, at
most, $\beta+\gamma;$ nevertheless, the test will consistently reject
the similarity model if it fails.
In the next section we show the performance in practice of this procedure.
Of course, one would like to control
\[
\limsup_n \sup_{(P,Q)\in H_0} \mathbb{P}_{(P,Q)}(p_{n,m}^*\leq\beta)
\]
instead of the bound given by our results. Some of the limitations of
our procedure come from
the smoothness requirements posed by our choice of metric, $\mathcal
{W}_2$. This could, perhaps,
be overcome with the use of the $L_1$ Wasserstein metric (but we would
lose the uniqueness
and consistency results given in Proposition \ref{uniqueness} and
Theorem \ref{consistencia})
and consideration of a less restrictive null hypothesis, $\tilde{H}_0$.
Uniformity in
$(P,Q)\in H_0$ is a more delicate issue, since one can take $P$ and $Q$
at an arbitrary (but positive)
Wasserstein distance from each other, but such that they are at
distance one in total variation.\vadjust{\goodbreak}
Perhaps a~different choice of metric could lead to some type of uniform
bound. We believe this
issue is worth further research.

Turning to the meaning of Condition (\ref{extrahipotesis}), observe
that the contaminations,
$P_1'$, $P_2'$, in the canonical decomposition in Proposition \ref
{equivalent} have disjoint
support but can be arbitrarily close in Wasserstein distance. With
Condition (\ref{extrahipotesis})
we avoid pathological cases in which some inconvenient distribution of
the contaminations allows
that some trimmings of $P$ and $Q$, with trimming size slightly above
the similarity level, are
close to each other without being too close to the common core.
Rather than pursuing an involved technical analysis
we include a couple of illustrative examples that show that the best possible
rate~$\rho$ depends on the \textit{degree of separation} between the
contaminating distributions $P_1'$,~$P_2'$ in the canonical
decomposition. In the well-separated case (when the distance between
the supports of $P_1'$ and $P_2'$ is positive), under additional
technical conditions we can take $\rho=1$ and we have that the optimal
trimming, $P_{n,\alpha_n}$, approaches
the common part, $P_0$, at the parametric rate: $\mathcal
{W}_2(P_{n,\alpha_n},P_0)=\mathrm{O}_P(n^{-1/2})$. Without this separation
we cannot take $\rho$ greater than $4/5$ and we have a nonparametric
rate of convergence:
$\mathcal{W}_2(P_{n,\alpha_n},P_0)=\mathrm{O}_P(n^{-2/5})$. Again, in our
examples we assume $P$ and $Q$ to have bounded support since
this is enough for applications.

%f3 ###
\begin{figure}

\includegraphics{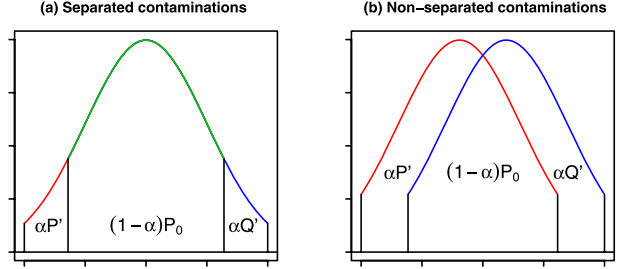}

\caption{Canonical decomposition in the separated (left) and
non-separated (right) cases.}
\label{separadasjuntas}
\end{figure}

\begin{Example}[(The well-separated case)]\label{casoseparado}
 Assume $P$ and $Q$ are probabilities on the real line with
quantile functions, $F^{-1}$ and $G^{-1}$, such that
$G^{-1}(t)=F^{-1}(t+\alpha)$, $0<t< 1-\alpha$ and $F^{-1}$ has a
bounded derivative (as in Figure \ref{separadasjuntas}\textup{(a)}). Then
$d_{\mathrm{TV}}(P,Q)=\alpha$ and, taking $\alpha_n=\alpha+\frac{K}{\sqrt{n}}$
for some $K>0$ and writing
$P_0$ for the common part in the canonical decomposition for $P$ and
$Q$, we have that
if $\tilde{P}_n\in\mathcal{R}_{\alpha_n}(P)$,
$\tilde{Q}_n\in\mathcal{R}_{\alpha_n}(Q),$ then
\[
\mathcal{W}_2(\tilde{P}_n,\tilde{Q}_n)=\mathrm{O}(n^{-1/2}) \quad \Rightarrow \quad \mathcal
{W}_2(\tilde{P}_n,P_0)=\mathrm{O}(n^{-1/2}).
\]
\end{Example}

\begin{Example}[(The non-separated case)]\label{casojunto}
We assume now that $P$ and $Q$ differ only in location and have a
symmetric, unimodal density. Without loss of generality, we write
$F(\cdot+\mu/2)$ and $F(\cdot-\mu/2)$
for the distribution functions of $P$ and $Q$, respectively, and $f$
for the density associated to~$F$. We suppose that $F$ has bounded
support and $f$ is strictly positive on it. Further, we assume $f$ to
be continuously differentiable with $f'<0$ in $(0,\sup(\operatorname{supp}(F)))$. If $\mu$ and $\alpha$ satisfy $1-\alpha=2(1-F(\mu
/2))=2F(-\mu/2),$ then $d_{\mathrm{TV}}(P,Q)=\alpha$ (see Figure \ref
{separadasjuntas}\textup{(b)}). If $\tilde{P}_n\in\mathcal{R}_{\alpha_n}(P)$,
$\tilde{Q}_n\in\mathcal{R}_{\alpha_n}(Q),$ then
\[
\mathcal{W}_2(\tilde{P}_n,\tilde{Q}_n)=\mathrm{O}(n^{-1/2}) \quad \Rightarrow \quad \mathcal
{W}_2(\tilde{P}_n,P_0)=\mathrm{O}(n^{-2/5}).
\]
\end{Example}

A proof of the claims in the last two examples is sketched in the \hyperref[appm]{Appendix}.

While this work is concerned mainly with testing $\alpha$-similarity in
two-sample problems, in many real problems
the interest could be focused on the estimation on the common core~$P_0$.
The results in Section \ref{SecPrelim}
ensure that the pooled probability, $P_{n,m}$, in our bootstrap
procedure is a consistent estimator of $P_0$ if
$\alpha$ equals the (unknown) distance in total variation between $P$
and $Q$. Our simulations in Section \ref{SecSimul}
(see Figure \ref{pvaloresoutliers} and the related comments) suggest
that the bootstrap $p$-value curves (the values of $p_{n,m}^*$ as a function
of~$\alpha$) change sharply from $0$ to $1$ around the true similarity
level. Maybe this rapid growth could
be used to give some estimation of the similarity level and, as a
result, of the common core. Further research is needed.

We conclude this section by presenting a simple upper bound for the
transportation cost between empirical measures.
This result, together with Theorem \ref{rate}, is the key in our proofs
of Theorems \ref{principal} and \ref{controldelnivel} and has some
independent interest. The proof is also included in the \hyperref[appm]{Appendix}. Here
$X_{1,1},\ldots,X_{1,n}; X_{2,1}, \ldots,X_{2,m}$ are i.i.d. $\mathbb
{R}^k$-valued random vectors with common distribution $P$ and
$Y_{1,1},\ldots,Y_{1,n}; Y_{2,1}, \ldots,Y_{2,m}$ are i.i.d. $Q$. We
write~$P_{n,1}$ and $P_{m,2}$ for the empirical measures based on
$X_{1,1},\ldots,X_{1,n}$ and $X_{2,1}, \ldots,X_{2,m}$, respectively,
and, similarly, $Q_{n,1}$ and $Q_{m,2}$ for the empirical measures
based on the $Y_{i,j}$.
Let us define
\[
S_{n,m}:=\mathcal{W}_p (P_{n,1},P_{m,2})\quad\mbox{and}\quad
T_{n,m}:=\mathcal{W}_p (Q_{n,1},Q_{m,2}).
\]
\begin{Prop}\label{empiricalcost}
With the above notation, if $p\geq1,$ then
\[
\mathcal{W}_p(\mathcal{L}(S_{n,m}),\mathcal{L}(T_{n,m}))\leq2 \mathcal
{W}_p (P,Q).
\]
\end{Prop}

%s4 ###
\section{Empirical analysis of the procedure}\label{SecSimul}
In this section we explore the performance of the procedure for finite
samples. The section is divided in two subsections that address the
analysis of a planned simulation study and of a case study,
respectively. To simplify our exposition we will assume equal sizes in
the two samples through the first subsection. All the computations have
been carried out with the programs available at
\url{http://www.eio.uva.es/\textasciitilde pedroc/R}.

%s4.1 ###
\subsection{A simulation study}\label{simulaciones}
We consider first an example that illustrates the over-fitting effect
on the bootstrap $p$-values. We generate 200 pairs of samples of size
$n=1000$ obtained from the $N(0,1)$ and the
$0.9N(0,1)+0.1N(10,3)$\vadjust{\goodbreak}
distributions. Then, for each pair of samples, we carry out the
bootstrap procedure (1000 bootstrap replicates in each run) for
trimming levels $\alpha=0.09$ and $0.11$. At this point an important
caution when dealing with mixtures should be made, namely the
distinction between the level (0.1 in our case) of the
``contaminating'' distribution in the mixture and the similarity level
between the non-contaminated and contaminated distributions. Of course,
both distributions are similar at level 0.1, but they are also similar
at a lower level (recall the canonical decomposition in Remark \ref
{desccanonica}). For example, since the supports of the $U(0,1)$ and
$U(1,2)$ distributions are disjoint, then the minimum level of
similarity between the $U(0,1)$ and $0.9U(0,1)+0.1U(1,2)$ distributions
is 0.1; but between the $N(0,1)$ and $0.9N(0,1)+0.1N(\mu,3)$
distributions, it is strictly lower for every $\mu$. For instance, this
level is 0.0484 if $\mu=0$, 0.0653 for $\mu=3$; or 0.0989 when $\mu=10$.

%f4 ###
\begin{figure}

\includegraphics{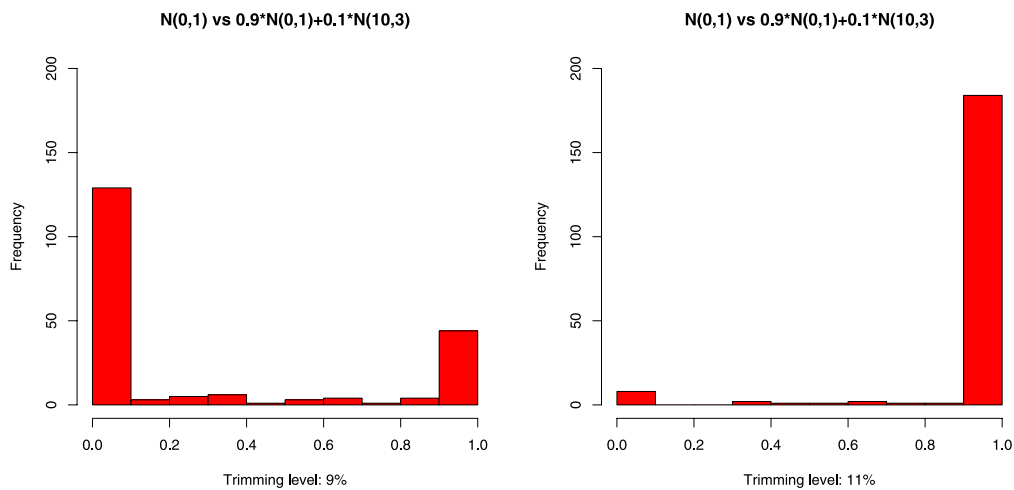}

\caption{Histograms, for different sizes of trimming, of the bootstrap
$p$-values obtained from 200 pairs of samples from $P=N(0,1)$ and
$Q=0.9N(0,1)+0.1N(10,3)$ distributions.}\vspace*{-1pt}
\label{grafico11}
\end{figure}

Figure \ref{grafico11} shows the absolute frequencies of the bootstrap
$p$-values, $p_{n,n}^*$, obtained in this example.

As stated above, the similarity level between the considered
distributions is 0.0989. Thus, the probability of obtaining an
observation from the non-common part in the mixture is 0.0989. Taking
into account sample sizes and the number of samples considered, the
expected number of times in which we obtain at most
110 `contaminating' observations in both samples is 158.13. In these
cases, after 0.11 trimming, we will be comparing similar samples and
should have no evidence against similarity. We note that 158 is
slightly below the observed frequency in the right bar of the right
histogram in Figure~\ref{grafico11}. On the other hand, the expected
number of times in which the amount of `contaminating' data exceeds 90
in both samples is 132.02. In this event, 0.09 trimming is unable to
remove contamination and
we should have strong evidence against similarity. We can check that
132 is close to the observed frequency in the left bar of the left
histogram in
Figure \ref{grafico11}.

\begin{figure}

\includegraphics{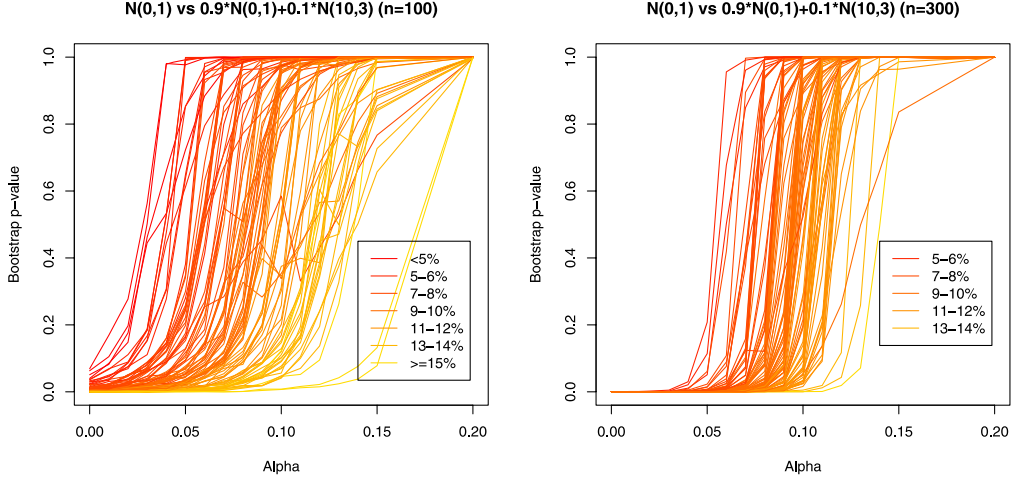}

\caption{Curves of bootstrap $p$-values obtained by varying the
trimming level ($\alpha$). Colors depend on the real proportion of data
coming from the $N(10,3)$ distribution in each particular sample. }
\label{pvaloresoutliers}
\end{figure}

The comments above suggest that the $p$-values are very sensitive to
the effective proportion of contamination in the data. This is further
illustrated with the plots in Figure~\ref{pvaloresoutliers}, which show
the curves of bootstrap $p$-values conditioned to different ranges of
contaminating proportion in the second sample (the amount of data
coming from the $N(10,3)$ distribution). In this figure we observe that
the transition from $p$-values close to 0 to $p$-values close to 1 is
very fast along the trimming level. In other words, the effect of
under-/over-trimming becomes apparent very quickly.

%f5 ###

We show next a simulation study to illustrate the power performance for
finite samples of the bootstrap procedure introduced in Section \ref
{bootstrapsection}, when the trimming level, $\alpha_n$, is determined
as in Theorem \ref{controldelnivel}. We consider two different cases,
comparing samples of the same size, $n$, of $P=N(0,1)$ versus $Q_i,
i=1,2$. In the first case, $Q_1=(1-\varepsilon)N(0,1)+\varepsilon
N(10,1)$; the contamination is due to outliers. In the second case, the
contamination is due to inliers and $Q_2=(1-\varepsilon
)N(0,1)+\varepsilon N(0,3)$. In both cases, the null hypothesis
is~$H_0\dvt\allowbreak d_{\mathrm{TV}}(P,Q_i)\le0.1$ and we use 1000 bootstrap pairs of samples
to obtain $p_{n,n}^*$, rejecting~$H_0$ if $p_{n,n}^* \le0.05 = \beta$.
Then we compute the rejection frequencies in 1000 iterations of the
procedure, obtaining the values shown in Tables \ref{potencia1} and \ref
{potencia2}. We do this for different values of $\varepsilon$ (then
different values of $\nu= d_{\mathrm{TV}}(P,Q_i)$) and different resampling
orders $n'=n^\rho$. The simulation shows that the bound given in
Theorem \ref{controldelnivel} is approached for moderate sizes in the
first case (see Table \ref{potencia1}, $\nu=0.10$). However, in the
second case, the procedure is conservative. The main conclusion is that
in both cases the contamination is detected, but detection is more
difficult in the case in which the contamination comes from inliers.

%t1 ###
\begin{table*}[t!]
\tabcolsep=0pt
\caption{Observed rejection frequencies for $H_0\dvt
d_{\mathrm{TV}}(P,Q_1)\le0.1$, $P=N(0,1)$, $Q_1=(1-\varepsilon
)N(0,1)+\varepsilon N(10,1)$,
where $\nu= d_{\mathrm{TV}}(P,Q_1)$ and $\beta=0.05$}\label{potencia1}
\begin{tabular*}{\textwidth}{@{\extracolsep{\fill}}ld{4.0}ld{1.3}d{1.3}d{1.3}d{1.3}d{1.3}d{1.3}d{1.3}d{1.3}@{}}
\hline
& &\multicolumn{1}{l}{{$\rho$:}} &\multicolumn{2}{l}{{1}}&\multicolumn{2}{l}{{4$/$5}}&\multicolumn{2}{l}{{2$/$3}}
&\multicolumn{2}{l@{}}{{1$/$2}} \\[-5pt]
& &\multicolumn{1}{l}{\hrulefill}&\multicolumn{2}{l}{\hrulefill}&\multicolumn{2}{l}{\hrulefill}&\multicolumn{2}{l}{\hrulefill}
&\multicolumn{2}{l@{}}{\hrulefill}\\
\multicolumn{1}{@{}l}{{$\nu$}}&\multicolumn{1}{l}{{$ n$}}
&\multicolumn{1}{l}{{$\gamma$:}}&\multicolumn{1}{l}{{0.05}} &\multicolumn{1}{l}{{0.01}}
&\multicolumn{1}{l}{{0.05}}&\multicolumn{1}{l}{{0.01}}&\multicolumn{1}{l}{{0.05}}&\multicolumn{1}{l}{{0.01}} &
\multicolumn{1}{l}{{0.05}} &\multicolumn{1}{l@{}}{{0.01}} \\
\hline
{0.10}&100 & &0.008 &0.001 &0.016 &0.003 &0.043 &0.006 &0.047 &0.007 \\
$\varepsilon\simeq0.10$&300 & &0.030 &0.007 &0.040 &0.015 &0.059 &0.017 &0.065 &0.019 \\
&1000& &0.052 &0.009 &0.092 &0.016 &0.098 &0.018 &0.114 &0.022 \\
[3pt]
{0.15}&100 & &0.130 &0.044 &0.207 &0.090 &0.246 &0.130 &0.252 &0.170 \\
$\varepsilon\simeq0.15$&300 & &0.587 &0.386 &0.648 &0.458 &0.687 &0.507 &0.703 &0.556 \\
&1000& &0.996 &0.980 &0.998 &0.985 &0.998 &0.986 &0.999 &0.990 \\
[3pt]
{0.20}&100 & &0.576 &0.403 &0.685 &0.515 &0.732 &0.585 &0.738 &0.624 \\
$\varepsilon\simeq0.20$&300 & &0.990 &0.973 &0.992 &0.981 &0.993 &0.985 &0.993 &0.986 \\
&1000& &1 &1 &1 &1 &1 &1 &1 &1 \\
[3pt]
{0.25}&100 & &0.919 &0.842 &0.953 &0.893 &0.969 &0.917 &0.970 &0.929 \\
$\varepsilon\simeq0.25$&300 & &1 &1 &1 &1 &1 &1 &1 &1 \\
&1000& &1 &1 &1 &1 &1 &1 &1 &1 \\
\hline
\end{tabular*}
\vspace*{-3pt}
\end{table*}

We close this subsection with a comparison to classical testing
procedures that could be adapted to the setup
of similarity testing. We recall from Proposition \ref{equivalent} that
testing $\alpha$-similarity of $P$ and $Q$ is equivalent to
testing whether $\sup_A |P(A)-Q(A)|\leq \alpha$, with~$A$ ranging
among all (measurable) sets. If we focus on sets
of type $A=(-\infty,x],$ then we could test the null hypothesis
$H_0\dvt
\sup_{x\in\mathbb{R}} |F(x)-G(x)|\leq \alpha$
using the Kolmogorov--Smirnov statistic: $D_n=\sup_{x\in\mathbb{R}}
|F_n(x)-G_n(x)|$, where $F_n$ and $G_n$ denote the
empirical distribution functions  (d.f.'s) based on the $X_i$ and the $Y_j$, respectively (and we
have assumed for simplicity samples of equal size). It is known
(see   \cite{Raghavachari1973}) that, provided $\sup_{x\in\mathbb{R}}
|F(x)-G(x)|=\lambda>0$, $\sqrt{n}(D_n-\lambda)$ converges weakly
to $Z_{\lambda}(F,G)=\max(Z_1,Z_2)$ with
\begin{eqnarray*}
Z_1&=&\sup_{\{x: F(x)-G(x)=\lambda\} } B_1\bigl(G(x)+\lambda\bigr)-B_2(G(x)),\\[-2pt]
Z_2&=&\sup_{\{x: G(x)-F(x)=\lambda\} } B_2(G(x))-B_1\bigl(G(x)-\lambda\bigr),
\end{eqnarray*}
where $B_1, B_2$ are independent Brownian bridges on $(0,1)$. With
standard arguments it can be shown that
$P(Z_{\lambda}(F,G)>t)\leq P(Z_{\lambda}>t)$ for $t>0$, with $Z_{\lambda
}=\sup_{0\leq x \leq1-\lambda} B_1(x+\lambda)-B_2(x)$.
Hence, if we choose $z_\alpha^{(\beta)}$ such that $P(Z_\alpha>z_\alpha
^{(\beta)})=\beta$, then the test that rejects when
\[
D_n>\alpha+\frac{1}{\sqrt{n}} z_\alpha^{(\beta)}
\]
is asymptotically of level $\beta$ for testing $H_0\dvt  \sup_{x\in\mathbb
{R}} |F(x)-G(x)|\leq \alpha$. The critical va\-lue~$z_\alpha^{(\beta)}$ can be approximated by Monte Carlo simulation. We
%
%t2 ###
\begin{table*}[t!]
\tabcolsep=0pt
\caption{Observed rejection frequencies for $H_0\dvt
d_{\mathrm{TV}}(P,Q_2)\le0.1$, $P=N(0,1)$, $Q_2=(1-\varepsilon)N(0,1) +
\varepsilon N(0,3)$, where $\nu= d_{\mathrm{TV}}(P,Q_2)$ and $\beta=0.05$}\label{potencia2}
\vspace*{-4pt}
\begin{tabular*}{\textwidth}{@{\extracolsep{\fill}}ld{4.0}ld{1.3}d{1.3}d{1.3}d{1.3}d{1.3}d{1.3}d{1.3}d{1.3}@{}}
\hline
& &\multicolumn{1}{l}{{$\rho$:}} &\multicolumn{2}{l}{{1}}&\multicolumn{2}{l}{{4$/$5}}
&\multicolumn{2}{l}{{2$/$3}}&\multicolumn{2}{l@{}}{{1$/$2}} \\[-6pt]
& &\multicolumn{1}{l}{\hrulefill}
&\multicolumn{2}{l}{\hrulefill}&\multicolumn{2}{l}{\hrulefill}&\multicolumn{2}{l}{\hrulefill}&\multicolumn{2}{l@{}}{\hrulefill}\\
\multicolumn{1}{@{}l}{{$\nu$}}&\multicolumn{1}{l}{{$ n$}}
&\multicolumn{1}{l}{{$\gamma$:}}&\multicolumn{1}{l}{{0.05}} &
\multicolumn{1}{l}{{0.01}}& \multicolumn{1}{l}{{0.05}}& \multicolumn{1}{l}{{0.01}}& \multicolumn{1}{l}{{0.05}}&\multicolumn{1}{l}{{0.01}} &
\multicolumn{1}{l}{{0.05}} & \multicolumn{1}{l@{}}{{0.01}} \\
\hline
 {0.10}&100 & &0 &0 &0 &0 &0 &0 &0 &0 \\
$\varepsilon\simeq0.21$&300 & &0 &0 &0 &0 &0 &0 &0 &0 \\
&1000& &0 &0 &0 &0 &0 &0 &0 &0 \\
[3pt]
 {0.15}&100 & &0.002 &0.000 &0.002 &0.001 &0.002 &0.001 &0.003 &0.001 \\
$\varepsilon\simeq0.31$&300 & &0.013 &0.003 &0.016 &0.005 &0.017 &0.006 &0.027 &0.008 \\
&1000& &0.185 &0.089 &0.196 &0.100 &0.210 &0.103 &0.235 &0.120 \\
[3pt]
{0.20}&100 & &0.037 &0.017 &0.048 &0.022 &0.060 &0.023 &0.065 &0.027 \\
$\varepsilon\simeq0.41$&300 & &0.397 &0.253 &0.418 &0.279 &0.437 &0.293 &0.490 &0.330 \\
&1000& &0.992 &0.979 &0.994 &0.979 &0.995 &0.982 &0.994 &0.983 \\
[3pt]
{0.25}&100 & &0.254 &0.146 &0.277 &0.163 &0.301 &0.189 &0.324 &0.195 \\
$\varepsilon\simeq0.52$&300 & &0.924 &0.846 &0.928 &0.856 &0.936 &0.866 &0.949 &0.888 \\
&1000& &1 &1 &1 &1 &1 &1 &1 &1 \\
[3pt]
{0.30}&100 & &0.565 &0.426 &0.599 &0.456 &0.629 &0.484 &0.654 &0.508 \\
$\varepsilon\simeq0.62$&300 & &0.996 &0.993 &0.998 &0.993 &0.998 &0.993 &0.999 &0.995 \\
&1000& &1 &1 &1 &1 &1 &1 &1 &1 \\
\hline
\end{tabular*}
\vspace*{-12pt}
\end{table*}
%t3 ###
\begin{table}
\tabcolsep=0pt
\tablewidth=290pt
\caption{Observed rejection frequencies for $H_0\dvt
d_{\mathrm{TV}}(P,Q)\le0.1$, $P=N(0,1)$, $Q=0.70 N(0,1)+0.15 N(2.35,1)+0.15
N(-2.35,1)$ at level $0.05$}\label{potencia3}
\vspace*{-4pt}
\begin{tabular*}{290pt}{@{\extracolsep{\fill}}lllll@{}}
\hline
$n$ & 100 & 300 & 500 & 1000\\
\hline
$D_n$ & 0.007 & 0.004 & 0.003 & 0.002\\
$\mathcal{W}_2$ & 0.007 & 0.091 & 0.320 & 0.875\\
\hline
\end{tabular*}
\vspace*{-4pt}
\end{table}
could try to use this procedure for testing the $\alpha$-similarity
model. Though, since we can find distributions that are arbitrarily
close in Kolmogorov--Smirnov distance but far from
each other in total variation distance, this alternative procedure can
fail badly. We show this in our last simulation study (see Table~\ref{potencia3}).
We have taken $P=N(0,1)$ and
$Q=0.70N(0,1)+0.15N(2.35,1)+0.15N(-2.35,1)$, a mixture with three
normal components. Here we have
$\sup_{x\in\mathbb{R}}|P(-\infty,x]-Q(-\infty,x]|=0.1$ and
$d_{\mathrm{TV}}(P,Q)=0.2$ and we test $H_0\dvt  d_{\mathrm{TV}}(P,Q)\leq0.1$ at level $0.05$.
We show the observed frequencies of rejection for $D_n$ and our
bootstrap procedure based on $\mathcal{W}_2$ as in Theorem~\ref{controldelnivel} with $\rho=4/5$,
$\gamma=0.01$. In this case we reject for bootstrap $p$-values larger
than $0.04$ to make the asymptotic probability of type I error less
than $0.05$.
We have considered sampling sizes $n=100, 300, 500$ and $1000$ and have
produced 10{,}000 replicates of the tests in each case.
We see that the Kolmogorov--Smirnov test fails to detect the
dissimilarity, even for large sample sizes, while the bootstrap
procedure suggested in this paper works reasonably for moderate
sizes.\vspace*{-3pt}

%s4.2 ###
\subsection{A case study}\label{markers}\vspace*{-3pt}

The data from this case study come from an admission exam to the
Universidad de Valladolid. 308 exams\vadjust{\goodbreak} on the same subject were randomly
assigned to 2 markers. The distribution of the exams was not exactly
balanced and markers received 152 and 156 exams, respectively. Each
exam was given a grade between 0 and 10 points. In the admission exams
some marking criteria are given to the markers with the goal of making
the grading process ``homogeneous''. The main goal of this study is to
determine whether the markers are using the same common criteria. Some
degree of deviation from this common pattern is allowed for each
marker. Therefore, we would like to assess the similarity of the
samples of marks for the different markers.

%f6 ###
\begin{figure}

\includegraphics{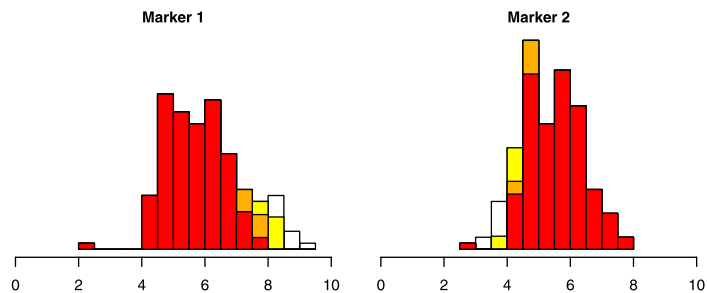}

\caption{Best trimmings between markers 1 and 2, in the example of
Section \protect\ref{markers}, $\alpha= 0.05$ (white), $\alpha= 0.10$
(white${}+{}$yellow) and $\alpha= 0.15$ (white${}+{}$yellow${}+{}$orange).}
\label{hist_correctores}
\end{figure}

%t4 ###
\begin{table*}
\tabcolsep=0pt
\caption{Bootstrap $p$-values arising from the
introduced bootstrap methodology,
applied to the similarity analysis between markers ($\beta=0.05$)}\label{pvaloresbootstrap}
\begin{tabular*}{\textwidth}{@{\extracolsep{\fill}}lld{1.3}d{1.3}d{1.3}d{1.3}d{1.3}d{1.3}d{1.3}d{1.3}@{}}
\hline
& \multicolumn{1}{l}{{$\rho$:}} &\multicolumn{2}{l}{{1}}&\multicolumn{2}{l}{{4$/$5}}
&\multicolumn{2}{l}{{2$/$3}}&\multicolumn{2}{l@{}}{{1$/$2}} \\[-5pt]
& \multicolumn{1}{l}{\hrulefill}
&\multicolumn{2}{l}{\hrulefill}&\multicolumn{2}{l}{\hrulefill}&\multicolumn{2}{l}{\hrulefill}&\multicolumn{2}{l@{}}{\hrulefill}\\
\multicolumn{1}{@{}l}{{$\alpha$}}
&\multicolumn{1}{l}{{$\gamma$:}}&\multicolumn{1}{l}{{0.05}}
&\multicolumn{1}{l}{{0.01}}&\multicolumn{1}{l}{{0.05}}&
\multicolumn{1}{l}{{0.01}}& \multicolumn{1}{l}{{0.05}}& \multicolumn{1}{l}{{0.01}} & \multicolumn{1}{l}{{0.05}} & \multicolumn{1}{l@{}}{{0.01}} \\
\hline
0 & &0 & 0 & 0 & 0 & 0 & 0 & 0 & 0 \\
0.05& &0.059 & 0.133 & 0.016 & 0.058 & 0.007 & 0.034 & 0.005 & 0.019 \\
0.10& &0.884 & 0.975 & 0.717 & 0.865 & 0.567 & 0.708 & 0.371 & 0.597 \\
0.15& &1 & 1 & 1 & 1 & 1 & 1 & 0.997 & 0.999 \\
0.20& &1 & 1 & 1 & 1 & 1 & 1 & 1 & 1 \\
\hline
\end{tabular*}
\end{table*}

The use of nonparametric methods strongly rejects, at level 0.05,
homogeneity between the considered marking distributions
(Wilcoxon--Mann--Whitney, $p$-value~$=
0.000$; and Kolmogo\-rov--Smirnov, $p$-value~$= 0.003$). In Figure \ref
{hist_correctores} we show the histograms corresponding to the full
data sets
and the progressive effects of best trimming, minimizing the Wasserstein
distance between the remaining subsample distributions. The white portions
of the bars represent the trimmed observations when the trimming size
is $\alpha= 0.05$, the union of the white and yellow portions are the
trimmed observations when $\alpha= 0.1$ and the orange portions
complete the trimming
corresponding to $\alpha= 0.15$. Notice that the best trimming is far
from being symmetric.

In Table \ref{pvaloresbootstrap} we have included the $p$-values
corresponding to the bootstrap
procedure introduced in Section \ref{bootstrapsection}. In every case,
for fixed $\beta=0.05$ and taking $\alpha_n$ as in Theorem~\ref
{controldelnivel}, we used 1000 bootstrap samples to compute the
$p$-values for the null hypothesis $H_0\dvt\allowbreak d_{\mathrm{TV}}(P,Q)\le\alpha$. In
general terms, these $p$-values show that both samples are not
$0.05$-similar, but they can be considered $0.10$-similar. The
considerations made in Section \ref{bootstrapsection} about Condition
(\ref{extrahipotesis}) show the convenience of using resampling orders
less than or equal to $n^{4/5}$, as we don't know if the supports of
the contaminating distributions are well separated or not.

\begin{appendix}
\setcounter{subsection}{0}
\section*{Appendix}\label{appm}
\renewcommand{\theequation}{\arabic{equation}}
\renewcommand{\thethm}{\arabic{thm}}
%s4.3 ###
\subsection{\texorpdfstring{Proof of Theorem \protect\ref{rate}}
{Proof of Theorem 2}}
Our proof is based on a parallel
result for the one-sample case. Let $P_n$ be the empirical measure
based on
i.i.d. random variables $X_1,\ldots,X_n$ with common distribution $P$.
In the particular case $P=Q$ and $\alpha=0$ we have $n\mathcal{W}^2_2(P_n
,Q)=\mathrm{O}_P(1)$ under sufficient integrability assumptions (see    \cite{Barrio}).
From the obvious bound $\mathcal{W}_2(\mathcal
R_\alpha(P_{n}) ,Q)\leq\mathcal{W}_2(P_n ,Q)$ we see that
$n\mathcal{W}^2_2(\mathcal R_\alpha(P_{n}) ,Q)=\mathrm{O}_P(1)$. Our first result
here shows that $n\mathcal{W}^2_2(\mathcal R_\alpha
(P_{n}),Q)=\mathrm{o}_P(1)$ even if $P\ne Q$.
\begin{thm}\label{unamuestra}
Assume that $Q \in\mathcal R_{\alpha_0} (P)$ for some $\alpha_0\in[0,1)$,
where $Q$ is supported in a~bounded interval, having a density function
that is bounded away from zero on its support, and with a bounded
derivative. If $\alpha_n \geq\alpha_0+r_n/\sqrt{n}$ for some sequence
$0\leq r_n\to\infty,$ then
\[
\sqrt{n}\mathcal{W}_2(\mathcal R_{\alpha_n} (P_{n}) ,Q) \to0   \qquad \mbox{in
probability as } n\to\infty.
\]
\end{thm}

\begin{pf} Arguing as in the proof of Proposition
\ref{equivalent} we can check that $Q\in
\mathcal R_{\alpha_0} (P)$ is equivalent to $P=(1-\alpha_0)Q+\alpha_0
P'$ for some
distribution $P'$. Hence, we can assume $X_n=(1-U_n)Y_n+U_nZ_n$, where
$\{Y_n\}_n$, $\{Z_n\}_n$ and $\{U_n\}_n$
are independent i.i.d. sequences with laws $Q$, $P'$ and Bernoulli with mean
$\alpha_0$, respectively. Write $N_n=\sum_{i=1}^n I(U_i=1)$. Then $N_n$
follows a binomial distribution with parameters $n$ and $\alpha_0$. Hence,
$\sqrt{n}(N_n/n-\alpha_0)\to\sqrt{\alpha_0(1-\alpha_0)}Z$, with $Z$
standard normal. We assume
w.l.o.g. that convergence holds, in fact, a.s. Write $n'=n-N_n$,
$\tilde{X}_1,\ldots,\tilde{X}_{n'}$
for the $Y_i$'s in the sample with associated $U_i=0$ (the uncontaminated
fraction of the sample:
$\tilde{X}_1,\ldots,\tilde{X}_{n'}$ are i.i.d. $Q$) and $\tilde
{P}_{n'}$ for
the empirical measure on the $\tilde{X}_i$'s. Observe that
$\tilde{P}_{n'}\in\mathcal{R}_{\tilde{\alpha}_n}(P_n)$ with
$\tilde{\alpha}_n=N_n/n$. Now we note that given $\alpha,\beta\in
[0,1)$, if
$Q \in\mathcal{R}_\alpha(P)$, then $\mathcal{R}_\beta( Q) \subset
{\mathcal R}_{\alpha+\beta- \alpha\beta} ( P)$. Hence,
$\mathcal{R}_{\hat{\alpha}_n}(\tilde{P}_{n'})\subset
\mathcal R_{\alpha_n} (P_{n})$ for
$\hat{\alpha}_n=(\alpha_n-\tilde{\alpha_n})/(\tilde{\alpha}_n)$ provided
$\alpha_n>\tilde{\alpha_n}$, which eventually holds. Consequently,
\[
\mathcal{W}_2 (\mathcal{R}_{\alpha_n}(P_n) , Q)\leq\mathcal{W}_2
(\mathcal{R}_{\hat{\alpha}_n}(\tilde{P}_{n'}) , Q).
\]
Thus, the result will follow if we prove it in the particular case
$P=Q$ and
$\alpha_0=0$.

We proceed in this case writing $F$ and $f$ for the distribution and density
functions of~$P$.
Recalling the parametrization in (\ref{caparameter}) we have
\[
\mathcal{W}^2_2(\mathcal{R}_{\alpha_n}(P_n) ,P)=\min_{h\in
\mathcal{C}_{\alpha_n}}\mathcal{W}^2_2((P_{n})_h ,P)=
\min_{h\in\mathcal{C}_{\alpha_n}} \int_0^1
\bigl(F_n^{-1}(h^{-1}(t))-F^{-1}(t)\bigr)^2\,\mathrm{d}t
\]
and we see that $n\mathcal{W}^2_2(\mathcal{R}_{\alpha_n}(P_n)
,P)=\min_{h\in\mathcal{C}_{\alpha_n}} M_n(h)$, where
\[
M_n(h)=\int_0^1
\biggl (\frac{\rho_n(t)}{f(F^{-1}(t))}-\sqrt
{n}\bigl(F^{-1}(h(t))-F^{-1}(t)\bigr) \biggr)
^2 h'(t)\,\mathrm{d}t
\]
and
$\rho_n(t)=\sqrt{n} f(F^{-1}(t)) (F_n^{-1}(t)-F^{-1}(t))$ is the weighted
quantile process.
Without loss of generality, we can assume that $\{X_n\}_n$ are defined
in a sufficiently rich
probability space in which there exist Brownian bridges, $B_n$, satisfying
%
%e13 ###
\begin{equation}\label{strongapproximation}
n^{1/2-\nu}\sup_{\fraca{1}{n}\leq t \leq
1-\fraca{1}{n}}\frac{|\rho_n(t)-B_n(t)|}{(t(1-t))^\nu}
= \cases{
 \mathrm{O}_P(\log n), &  \quad $\mbox{if } \nu=0$,\cr
\mathrm{O}_P(1), &  \quad $\mbox{if } 0<\nu\leq1/2$
 }
\end{equation}
(this is guaranteed by Theorem 6.2.1 in  \cite{CsH93}). Now, defining
\[
\tilde{N}_n(h)=\int_0^1
 \biggl(\frac{B_n(t)}{f(F^{-1}(t))}-\sqrt{n}\bigl(F^{-1}(h(t))-F^{-1}(t)\bigr) \biggr)^2
h'(t)\,\mathrm{d}t,
\]
and assuming w.l.o.g. that $\alpha_n\leq1-\delta$ for some $\delta>0$ we
have that
\[
\sup_{h_\in\mathcal{C}_\alpha} |M_n(h)^{1/2}- \tilde{N}_n(h)^{1/2}|\leq
 \biggl(\frac{1}{\delta}\int_0^1
 \biggl(\frac{\rho_n(t)-B_n(t)}{f(F^{-1}(t))} \biggr)^2\,\mathrm{d}t
 \biggr)^{1/2}=\mathrm{o}_P(1).
\]
The last equality follows from (\ref{strongapproximation}), taking $\nu=0$,
because, since $f$ is bounded below
\[
\int_{1/n}^{1-1/n}  \biggl(\frac{\rho_n(t)-B_n(t)}{f(F^{-1}(t))} \biggr)
^2\,\mathrm{d}t\leq\frac{\log n}{\sqrt{n}} \int_0^1\frac{1}{f^2(F^{-1}(t))}\,\mathrm{d}t \mathrm{O}_P(1)
=\mathrm{o}_P(1).
\]
Thus, the conclusion will follow if we show
$\min_{h\in\mathcal{C}_{\alpha_n}} \tilde{N}_n(h) \to0$ in probability or,
equivalently,
if we show that $\min_{h\in\mathcal{C}_{\alpha_n}} {N}_n(h) \to0$ in
probability, where
\[
{N}_n(h)=\int_0^1
 \biggl(\frac{B(t)}{f(F^{-1}(t))}-\sqrt{n}\bigl(F^{-1}(h(t))-F^{-1}(t)\bigr) \biggr)^2
h'(t)\,\mathrm{d}t
\]
and $B$ is a fixed Brownian bridge. To check that
$\min_{h\in\mathcal{C}_{\alpha_n}} {N}_n(h) \to0$ in probability, we observe
that $\min_{h\in\mathcal{C}_{\alpha_n}} {N}_n(h)\leq\frac{1}{\delta}
\min_{k\in\mathcal{G}_{n}} R_n(k)$, where
\begin{eqnarray*}
R_n(k)=\int_0^1
\biggl (\frac{B(t)}{f(F^{-1}(t))}-\sqrt{n}\bigl(F^{-1}\bigl(t+k(t)/\sqrt{n}\bigr)-F^{-1}(t)\bigr)
 \biggr)^2\,\mathrm{d}t
\end{eqnarray*}
and $\mathcal{G}_{n}$ is the set of real-valued, absolutely continuous
functions
on $[0,1]$ such that $k(0)=k(1)=0$ and $-\sqrt{n}\leq k'(t)\leq r_n$ for
almost every $t$. We assume w.l.o.g. $r_n\leq r_{n+1}$ for every~$n$.
Then $\mathcal{G}_{n}\subset\mathcal{G}_{n+1}$ for every $n$
and $\mathcal{G}:=\bigcup_{n\geq1}\mathcal{G}_{n}$ is the
set of all absolutely continuous\vspace*{1pt} functions
on $[0,1]$ such that $k(0)=k(1)=0$ and $k'$ is (essentially) bounded. From
our hypotheses it follows easily that,
for $k\in\mathcal{G}$,
\[
R_n(k)\to R(k):=
\int_0^1  \biggl(\frac{B(t)-k(t)}{f(F^ {-1}(t))} \biggr)^2\,\mathrm{d}t
\]
and hence $\min_{k\in\mathcal{G}_{n}} {R}_n(k) \to0$ (therefore
$n\mathcal{W}_2^2(\mathcal{R}_{\alpha_n}(P_n),P)\to0$) will follow if we
show that $\inf_{k\in\mathcal{G}} R(k)=0$. But this can be checked easily
by noting, for instance, that if $k_n$ is the function that
interpolates $B(t)$ at knots $i/n$, $i=0,\ldots,n$, and is linear in
between, then
we have $k_n\in\mathcal{G}$ and $R(k_n)\to0$.
\end{pf}

\begin{pf*}{Proof of Theorem \ref{rate}}
We write $\alpha_0=d_{\mathrm{TV}}(P,Q)$ and take $P_0$ as in the canonical decomposition
in Proposition \ref{equivalent} (we take $\mu$ to be the Lebesgue
measure there). Then
$P_0 \in\mathcal R_{\alpha_0} (P)$ holds with $P$ and $P_0$ playing the
roles of $P$ and $Q$ and the density of~$P_0$ satisfies the
assumptions in Theorem~\ref{unamuestra} (in fact $f_0=(f\wedge
g)/(1-\alpha_0)$ has a bounded derivative a.e., but this suffices for the
strong approximation in the proof of Theorem~\ref{unamuestra}). Hence,
$\sqrt{n}\mathcal{W}_2(\mathcal{R}_{\alpha_n}(P_n), P_0)\to0$ in
probability and similarly for
$\sqrt{n}\mathcal{W}_2(\mathcal{R}_{\alpha_n}(Q_n), P_0)$. The triangle
inequality for $\mathcal{W}_2$ yields the conclusion.
\end{pf*}

%s4.4 ###
\subsection{Asymptotic theory for the bootstrap}
The behavior of the bootstrap $p$-value under the alternative follows
from the next result.

\begin{Prop}\label{controlalternativa}
Assume $X_{n,1},\ldots,X_{n,n'}; Y_{n,1}, \ldots,Y_{n,m'}$ are i.i.d.
random variables with common distribution $P_n\in\mathcal{F}_2$
such that $\mathcal{W}_2(P_n,P)\to0$. If $P_{n'}^*$ and $Q_{m'}^*$
denote the empirical measures on
$X_{n,1},\ldots,X_{n,n'}$ and $Y_{n,1}, \ldots,Y_{n,m'}$,
respectively, and $n', m'\to\infty,$ then
\[
\mathcal{W}_2(P_{n'}^*,Q_{m'}^*)\to0 \qquad\mbox{in probability}.
\]
\end{Prop}

\begin{pf} By Proposition \ref{empiricalcost} it is
enough to consider the case $P_n=P$ for all $n$. But then
$P_{n'}\to_w P$ a.s. by the Glivenko--Cantelli theorem while the law of
large numbers gives convergence of second-order moments. These two
facts imply that
$\mathcal{W}_2(P_{n'}^*,P)\to0$ (and for $\mathcal{W}_2(Q_{m'}^*,P)$
as well).
\end{pf}

Now we take care of the null hypothesis. The next result will be useful
for $P$ and $Q$ away from the boundary. Its proof is analogous to that
of Theorem 2.1 in   \cite{Bickel}.
\begin{Prop}\label{BF}
Assume $X_{n,1},\ldots,X_{n,n'}$ are i.i.d. random variables with
common distribution $P_n\in\mathcal{F}_2$
such that $\mathcal{W}_2(P_n,P)\to0$. If $\bar{X}_{n,n'}:=\frac{1}{n'}
\sum_{i=1}^{n'} X_{n,i}$, then
\[
\sqrt{n'}(\bar{X}_{n,n'}- \mu_n) \to_w N(0,\sigma^2),
\]
where $\mu_n = E(\bar{X}_{n,n'})$ and $\sigma^2$ is the variance of $P$.
\end{Prop}

\begin{pf*}{Proof of Theorem \ref{principal}} We will assume for
simplicity $n=m$ and $n'=m'$. The general case can be
handled with straightforward modifications. We consider first the case
$d_{\mathrm{TV}}(P,Q)>\alpha$. In this case we have
(Theorem \ref{consistencia}) that $\mathcal{W}_2(P_{n,\alpha_n},P_\alpha
)\to0$ and $\mathcal{W}_2(Q_{n,\alpha_n},Q_\alpha)\to0$ a.s. Since
\[
\mathcal{W}_2^2\bigl(aP_1+(1-a)P_2,aQ_1+(1-a)Q_2\bigr)\leq a\mathcal
{W}_2^2(P_1,Q_1)+(1-a)\mathcal{W}_2^2(P_2,Q_2)
\]
for probabilities $P_i,Q_i\in\mathcal{F}_2$ and $a\in[0,1]$
(see    \cite{Pedro2}) it follows that $\mathcal{W}_2(R_{n,n},\lambda P_\alpha
+(1-\lambda)Q_\alpha)\to0$ a.s. Note that
\[
p_{n,n}^*=\mathbb{P}^* \Biggl(\mathcal{W}_2(P_{n'}^*,Q_{n'}^*)>\sqrt{\frac
{n}{n'}}\mathcal{W}_2(P_{n,\alpha_n},Q_{n,\alpha_n}) \Biggr).
\]
Now, Theorem \ref{consistencia} implies that $\mathcal{W}_2(P_{n,\alpha
_n},Q_{n,\alpha_n})\to\mathcal{W}_2(\mathcal{R}_\alpha(P),\mathcal
{R}_\alpha(Q))>0$,
while $n/n'$ is bounded away from 0 by assumption. This, together with
Proposition \ref{controlalternativa}, gives (ii).

We assume now that $d_{\mathrm{TV}}(P,Q)<\alpha$. Then Theorem \ref{rate}
ensures that $\sqrt{n} \mathcal{W}_2(P_{n,\alpha_n},\allowbreak Q_{n,\alpha
_n})\to0$ in probability.
Now, if $P_1,P_2$ are probabilities in $\mathcal{F}_2$ with means $\mu
_1,\mu_2$ and~$\bar{P}_1,\bar{P}_2$ are their centered versions, then
it is easy to check that
$\mathcal{W}_2^2(P_1,P_2)=(\mu_1-\mu_2)^2+\mathcal{W}_2^2(\bar{P}_1,\bar
{P}_2)$ and, therefore, $\mathcal{W}_2^2(P_1,P_2)\geq(\mu_1-\mu_2)^2$.
Let $\bar{X}_{n'}^*$ and $\bar{Y}_{n'}^*$, respectively, denote the
means corresponding to the $X$'s and $Y$'s bootstrap samples, and~$\mu
_{n} $ be the mean of the parent bootstrap distribution, $R_{n,n}$. Then
\[
n'\mathcal W_2^2 (P_n^{*},Q_m^{*} )\geq n' (\bar
{X}_{n'}^*-\bar{Y}_{n'}^* )^2= \bigl(\sqrt{n'}(\bar{X}_{n'}^*-\mu
_{n} )-\sqrt{n'}(\bar{Y}_{n'}^*-\mu_{n} ) \bigr)^2.
\]
%
%Therefore, by conditional independence of these means, (i) will hold
%if we prove that any subsequence of $\sqrt n(\bar{X}_n^*-\mu_{n,m} )$
%and of $\sqrt%n(\bar{Y}_m^*-\mu_{n,m})$ have convergent subsequences
%with normal limit laws (possibly degenerated) whose sum of variances
%is at least $\sigma_0^2$.

From the Glivenko--Cantelli theorem we have a.s. tightness of $\{P_{n}\}
_{n}$ and $\{Q_{n}\}_{n}$ and, as a consequence, of $P_{n,\alpha_n}$
and $Q_{n,\alpha_n}$ (see Proposition 2.1 in  \cite{Pedro2}). We can
assume, taking subsequences if necessary, that $P_{n,\alpha_n}\to_w
P_0$ and
$Q_{n,\alpha_n}\to_w Q_0$ for some probabilities~$P_0,\allowbreak Q_0$. A little
thought shows that, necessarily, $P_0\in\mathcal{R}_\alpha(P)$ and
$Q_0\in\mathcal{R}_\alpha(Q)$. Since $\mathcal{W}_2(P_{n,\alpha
_n},Q_{n,\alpha_n}) \to0$, necessarily, $P_0=Q_0 \in\mathcal{R}_\alpha
(P)\cap\mathcal{R}_\alpha(Q)$. Also, since $P,Q\in\mathcal F_2$, the
strong law of large numbers shows that the map $x^2$ is uniformly
integrable with respect to $\{P_n\}_n$ and $\{Q_n\}_n$ a.s., hence also
with respect to $\{P_{n,\alpha_n}\}_{n}$ and $\{Q_{n,\alpha_n}\}_{m}$.
Thus, perhaps through subsequences, $\mathcal{W}_2(P_{n,\alpha
_n},P_0)\to0$ and $\mathcal{W}_2(Q_{n,\alpha_n},P_0)\to0$, hence
$\mathcal{W}_2(R_{n,n},P_0)\to0$ for some $P_0\in\mathcal{R}_\alpha
(P)\cap\mathcal{R}_\alpha(Q)$.

The function that sends $P$ to its variance is continuous in $\mathcal
{F}_2$ for the $\mathcal{W}_2$ metric. Hence, since $\mathcal{R}_\alpha
(P)\cap\mathcal{R}_\alpha(Q)$ is compact, the variance attains its
minimum there. Let us write
$\sigma_0^2=\min_ {R\in\mathcal{R}_\alpha(P)\cap\mathcal{R}_\alpha
(Q)} \operatorname{{Var}} (R)$. Then $\sigma_0>0$ (a trimming of a probability
with a density has a density, hence, cannot have null variance) and if
we write $\sigma^2$ for the variance of~$P_0$, we have
\begin{eqnarray*}
p_{n,n}^*&=&\mathbb{P}^* \bigl(\sqrt{n'}\mathcal{W}_2
(P_{n'}^*,Q_{n'}^*)> \sqrt{n} \mathcal{W}_2(P_{n,\alpha_n},Q_{n,\alpha
_n})  \bigr)\\
&\geq& \mathbb{P}^* \biggl(\biggl|{{\frac{\sqrt{n'}}{2\sigma}} }(\bar
{X}_{n'}^*-\bar{Y}_{n'}^*)\biggr| > {{\frac{\sqrt{n}}{2\sigma}}}
\mathcal{W}_2(P_{n,\alpha_n},Q_{n,\alpha_n})  \biggr)\\
&\geq& \mathbb{P}^* \biggl(\biggl|{{\frac{\sqrt{n'}}{2\sigma}} }(\bar
{X}_{n'}^*-\bar{Y}_{n'}^*)\biggr| > {{\frac{\sqrt{n}}{2\sigma_0}}}
\mathcal{W}_2(P_{n,\alpha_n},Q_{n,\alpha_n})  \biggr).
\end{eqnarray*}
Thus, Proposition \ref{BF} and the fact that $\sqrt{n} \mathcal
{W}_2(P_{n,\alpha_n},Q_{n,\alpha_n})\to0$ yield that $p_{n,n}^*\to1$
in probability, showing (i).
\end{pf*}

\begin{pf*}{Proof of Theorem \ref{controldelnivel}}
As in the proof of Theorem \ref{rate}, we assume that
$X_n=(1-U_n)A_n+U_n B_n$, $Y_n=(1-V_n)C_n +V_n D_n$ with
$\{A_n\}_n$, $\{B_n\}_n$, $\{C_n\}_n$, $\{D_n\}_n$, $\{U_n\}_n$, $\{
V_n\}_n$ independent i.i.d. sequences of which
$\{A_n\}_n$ and $\{C_n\}_n$ have common distribution $P_0$ while $\{
U_n\}_n$ and $\{V_n\}_n$ are Bernoulli with mean $\alpha$. We write
$N_n=\sum_{i=1}^n I(U_i=1)$ and
$M_n=\sum_{i=1}^n I(V_i=1)$. Also we put $n_1'=n-N_n$, $n_2'=n-M_n$ and
write $\tilde{X}_1,\ldots,\tilde{X}_{n_1'}$ and $\tilde{Y}_1,\ldots
,\tilde{Y}_{n_2'}$ for the data corresponding to $U_i=0$ and $V_i=0$,
respectively.

On the set $E_n:=(N_n\leq n\alpha_n, M_n\leq n\alpha_n)$, the empirical
measures on $\tilde{X}_1,\ldots,\tilde{X}_{n_1'}$ and $\tilde
{Y}_1,\ldots,\tilde{Y}_{n_1'}$ (which we denote $\tilde{P}_{n_1'}$ and
$\tilde{Q}_{n_2'}$) satisfy $\tilde{P}_{n_1'}\in\mathcal{R}_{\alpha
_n}(P_n)$ and $\tilde{Q}_{n_2'}\in\mathcal{R}_{\alpha_n}(Q_n)$. Hence,
we have $\mathcal{W}_2(P_{n,\alpha_n},Q_{n,\alpha_n})\leq\mathcal{W}_2
(\tilde{P}_{n_1'},\tilde{Q}_{n_2'})$. Thus,
\begin{eqnarray*}
\mathbb{P}(p_{n,n}^*\leq\beta)\leq \mathbb{P}(E_n^C)+ \mathbb
{P}\bigl((\tilde{p}_n^*\leq\beta)\cap E_n\bigr),
\end{eqnarray*}
where
\[
\tilde{p}_n^*=\mathbb{P}^*  \bigl( \sqrt{n'} \mathcal
{W}_2(P_{n'}^{*},Q_{n'}^{*}) > \sqrt{n(1-\alpha)} \mathcal{W}_2 (\tilde
{P}_{n_1'},\tilde{Q}_{n_2'}) \bigr).
\]
By the central limit theorem (CLT) we have $\mathbb{P}(E_n^C)\to\gamma
$. Hence it suffices to control $\mathbb{P}((\tilde{p}_n^*\leq\beta
)\cap E_n)$. If $J_1,\ldots,J_{n'}$,
$L_1,\ldots,L_{n'}$ are i.i.d. random variables with law $P_0$,
independent of the data (both original and bootstrap) and $\mu_{n'}$,
$\nu_{n'}$ are the empirical measures, then Theorem \ref{empiricalcost}
and the fact that $\mathcal{W}_2(\mathcal{L}(aX),\mathcal
{L}(aY))=a\mathcal{W}_2(\mathcal{L}(X),\mathcal{L}(Y))$ for $a>0$ imply
\[
\mathcal{W}_2\bigl(\mathcal{L}^*\bigl(\sqrt{n'} \mathcal
{W}_2(P_{n'}^{*},Q_{n'}^{*})\bigr) , \mathcal{L}\bigl(\sqrt{n'} \mathcal{W}_2(\mu
_{n'},\nu_{n'})\bigr)\bigr)\leq2\sqrt{n'}\mathcal{W}_2(R_{n,n},P_0).
\]
By Lemma \ref{lematecnico} below $\sqrt{n'}\mathcal
{W}_2(R_{n,n},P_0)I_{E_n}\to0$ in probability. The assumptions on $P$
and~$Q$ yield that
$\sqrt{n'} \mathcal{W}_2(\mu_{n'},\nu_{n'})$ converges weakly to a
non-null limiting distribution as in~(\ref{weaklimit})
(with a proof as in Theorem 4.6 in   \cite{Barrio}). We
call $\eta$ the limit probability measure. Then
\[
\bigl|\tilde{p}_n^*-\eta\bigl(\bigl(\sqrt{n(1-\alpha)} \mathcal{W}_2 (\tilde
{P}_{n_1'},\tilde{Q}_{n_2'}),\infty\bigr)\bigr)\bigr|I_{E_n}\to0
\]
in probability. As a consequence,
\[
\mathbb{P}\bigl((\tilde{p}_n^*\leq\beta)\cap E_n\bigr)-\mathbb{P} \bigl(\bigl(\eta
\bigl(\bigl(\sqrt{n(1-\alpha)} \mathcal{W}_2 (\tilde{P}_{n_1'},\tilde
{Q}_{n_2'}),\infty\bigr)\bigr)\leq\beta\bigr)\cap E_n \bigr)\to0.
\]
But\vspace*{-1pt}
\begin{eqnarray*}
&& \mathbb{P} \bigl(\bigl(\eta\bigl(\bigl(\sqrt{n(1-\alpha)} \mathcal{W}_2 (\tilde
{P}_{n_1'},\tilde{Q}_{n_2'}),\infty\bigr)\bigr)\leq\beta\bigr)\cap E_n \bigr)\\
&& \quad \leq
\mathbb{P} \bigl(\bigl(\eta\bigl(\bigl(\sqrt{n(1-\alpha)} \mathcal{W}_2 (\tilde
{P}_{n_1'},\tilde{Q}_{n_2'}),\infty\bigr)\bigr)\leq\beta\bigr) \bigr)\to\beta,
\end{eqnarray*}
since, as above, $\sqrt{n(1-\alpha)} \mathcal{W}_2 (\tilde
{P}_{n_1'},\tilde{Q}_{n_2'})$ converges weakly to $\eta$. This
completes the proof.
\end{pf*}

The following technical result has been used in the proof of Theorem
\ref{controldelnivel}.
\begin{Lemm}\label{lematecnico}
With the notation and assumptions of Theorem \ref{controldelnivel},
\[
\sqrt{n'}\mathcal{W}_2(R_{n,n},P_0)I_{E_n}=\mathrm{o}_P(1).
\]
\end{Lemm}

\begin{pf} We use the parametrization in (\ref
{caparameter}). We have $P_{n,\alpha_n}=(P_{n})_{h_n}$, $Q_{n,\alpha
_n}=(Q_{n})_{l_n}$,
for some $h_n, l_n\in\mathcal{C}_{\alpha_n}$. Writing $F_{n}^{-1}$,
$G_{n}^{-1}$, $F^{-1}$ and $G^{-1}$ for the quantile functions of
$P_n$, $Q_n$, $P$ and $Q$
we have $\mathcal{W}^2(P_{n,\alpha_n},Q_{n,\alpha_n})=\|F_n^{-1}\circ
h_n^{-1}-G_n^{-1}\circ l_n^{-1}\|_2$, with $\|\cdot\|_2$ denoting the
usual norm in $L_2 (0,1)$, namely,
$\|b\|_2^2=\int_0^1 b^2$. Now
\begin{eqnarray*}
 &&\| (F_n^{-1}\circ h_n^{-1}-G_n^{-1}\circ
l_n^{-1})-(F^{-1}\circ h_n^{-1}-G^{-1}\circ l_n^{-1})\|_2\\
&& \quad \leq\| F_n^{-1}\circ h_n^{-1}-F^{-1}\circ h_n^{-1}\|_2+\|
G_n^{-1}\circ l_n^{-1}-G^{-1}\circ l_n^{-1}\|_2\\
&& \quad \leq\frac1{\sqrt{1-\alpha_n}}( \| F_n^{-1}-F^{-1}\|_2+\|
G_n^{-1}-G^{-1}\|_2),
\end{eqnarray*}
where we have used that $\int_0^1
(F^{-1}(h^{-1}(t))\!-\!G^{-1}(h^{-1}(t))^2\,\mathrm{d}t\!=\!\int_0^1
(F^{-1}(x)\!-\!G^{-1}(x)^2 h'(x)\,\mathrm{d}x$. The assumptions on $P$ and $Q$ ensure that,
as in (\ref{weaklimit}), $ \| F_n^{-1}-F^{-1}\|_2+\|G_n^{-1}-G^{-1}\|
_2=\mathrm{O}_P(n^{-1/2})$. On the other hand, on $E_n$,
\[
\| F_n^{-1}\circ h_n^{-1}-G_n^{-1}\circ l_n^{-1}\|_2=\mathcal
{W}_2(P_{n,\alpha_n},Q_{n,\alpha_n})\leq\mathcal{W}_2 (\tilde
{P}_{n_1'},\tilde{Q}_{n_2'})=\mathrm{O}_P(n^{-1/2}).
\]
Combining these two facts we see that
$\mathcal{W}_2(P_{h_n},Q_{h_n})I_{E_n}=\|F^{-1}\circ
h_n^{-1}-G^{-1}\circ l_n^{-1}\|_2 I_{E_n}=\mathrm{O}_P(n^{-1/2}).$ Using (\ref
{extrahipotesis}) we see that $\mathcal{W}_2(P_{h_n},P_0)=\mathrm{O}(n^{-\rho
/2})$. Since $\mathcal{W}_2(P_{h_n},P_{n,\alpha_n})= \mathrm{O}_P(n^{-1/2})$, we
conclude that
$\mathcal{W}_2(P_{n,\alpha_n},P_0)I_{E_n}=\mathrm{O}(n^{-\rho/2})$. Convexity
and a similar argument for $Q_{n,\alpha_n}$ yield the result.
\end{pf}

\begin{pf*}{Proof of Example \ref{casoseparado}}
The fact that $d_{\mathrm{TV}}(P,Q)=\alpha$ follows from noting (with some abuse
of notation) that for $\tilde{F}^{-1}\in\mathcal{R}_\alpha(P)$ and
$\tilde{G}^{-1}\in\mathcal{R}_\alpha(Q)$
\[
\tilde{F}^{-1}(t)\leq F^{-1}\bigl(\alpha+(1-\alpha)t\bigr)\leq\tilde{G}^{-1}(t).
\]
Hence, the probability $P_0$ with quantile $F_0^{-1}(t)=F^{-1}(\alpha
+(1-\alpha)t)$ is the unique element in $\mathcal{R}_\alpha(P)\cap
\mathcal{R}_\alpha(Q)$.
Next we observe that, for $\tilde{F}^{-1}\in\mathcal{R}_{\alpha_n}(P)$,
\begin{eqnarray*}
F^{-1}(t)&\leq& F^{-1}\bigl(\alpha_n+(1-\alpha_n)t\bigr) \\
&\leq& F_0^{-1}(t)+\bigl(F^{-1}\bigl(\alpha_n+(1-\alpha_n)t\bigr)-F^{-1}\bigl(\alpha
+(1-\alpha_n)t\bigr)\bigr).
\end{eqnarray*}
Similarly, if $\tilde{G}^{-1}\in\mathcal{R}_{\alpha_n}(Q)$, $\tilde
{G}^{-1}(t)\geq F_0^{-1}(t)-(F^{-1}(\alpha_n+(1-\alpha
_n)t)-F^{-1}(\alpha+(1-\alpha_n)t))$
and, combining both inequalities, we get $|F_0^{-1}(t)-\tilde
{F}^{-1}(t)|\leq|\tilde{F}^{-1}(t)-\tilde{G}^{-1}(t)|+|F^{-1}(\alpha
_n+(1-\alpha_n)t)-F^{-1}(\alpha+(1-\alpha_n)t)|$
and the bound follows from the triangle inequality.
\end{pf*}

\begin{pf*}{Proof of Example \ref{casojunto}}
We write
$F_0$ for the distribution function of $P_0$, hence, $F_0^{-1}(y)=\mu
/2+F^{-1} ((1-\alpha)y)$ for $y\in(0,1/2]$ and $F_0^{-1}(y)=-\mu
/2+F^{-1} (\alpha+(1-\alpha)y)$ for
$y\in[1/2,1)$. Similarly, we write $\tilde{F}_n$ and $\tilde{G}_n$
for
the distribution functions of $\tilde{P}_n$ and $\tilde{Q}_n$,\vspace*{1pt} respectively.\vadjust{\goodbreak}
Necessarily, $\tilde{P}_n(0,\infty)\leq\frac1 {1-\alpha_n}(1-F(\frac
\mu2))=\frac1 2  (1+ \frac K {(1-\alpha_n)\sqrt{n}} )$. We
write $\beta_n={\frac1 2}-\tilde{P}_n(0,\infty)$.
It follows from the fact that $\mathcal{W}_2(\tilde{P}_n,\tilde
{Q}_n)\to0$ that $\mathcal{W}_2(\tilde{P}_n,P_0)\to0$ and, therefore,
that $\beta_n\to0$.
We give next a lower bound for $\mathcal{W}_2(\tilde{P}_n,\tilde
{Q}_n),$ assuming that $\beta_n>0$.
If this is the case
%
%e14 ###
\begin{equation}\label{ineq4}
\tilde{F}_n^{-1}(t)\leq-\frac\mu2 +F^{-1}\biggl(\alpha+(1-\alpha_n)(t-\beta
_n)+{\frac K {2\sqrt{n}}}\biggr), \qquad t\in\biggl(0,{\frac
1 2}+\beta_n\biggr).
\end{equation}
On the other hand $\tilde{G}_n^{-1}((1-\alpha_n)t)\geq\mu/2
+F^{-1}((1-\alpha_n)t)$. Standard computations show that there is a
unique $a=a(\beta_n)>0$ such that $F(a-\frac{\mu}2)-F(a+\frac{\mu
}2)+\alpha=(1-\alpha)\beta_n$ and that
\[
-\frac\mu2 +F^{-1}\bigl(\alpha+(1-\alpha)(t-\beta)\bigr)\leq\mu/2
+F^{-1}\bigl((1-\alpha)t\bigr)
\]
for $t\in(\frac1{1-\alpha} F(-a-\frac \mu2 ),\frac1 2)$. From this
we get that
%
%e15 ###
\begin{equation}\label{lalowerbound}
\mathcal{W}_2(\tilde{P}_n,\tilde{Q}_n)\geq\sqrt{g_1(\beta_n)}-s_{n,1}-s_{n,2},
\end{equation}
where $g_1(\beta)=\int_{F(-a-\mu/2)/(1-\alpha)}^{1/2} (\mu
+F^{-1}((1-\alpha)t)-F^{-1}(\alpha+(1-\alpha)(t-\beta)))^2\,\mathrm{d}t$,
$s_{n,1}^2=\int_{F(-a-\mu/2)/(1-\alpha)}^{1/2} (F^{-1}((1-\alpha
)t)-F^{-1}((1-\alpha_n)t))^2\,\mathrm{d}t$,
$s_{n,2}^2=\int_{F(-a-\mu/2)/(1-\alpha)}^{1/2} (F^{-1}(\alpha+(1-\alpha
)(t-\beta_n))-F^{-1}(\alpha+(1-\alpha_n)(t-\beta_n)+{\frac K
{2\sqrt{n}}}))^2\,\mathrm{d}t.$
A routine use of Taylor expansions yields $\lim_{\beta\to0+}\frac
{g_1(\beta)}{\beta^{5/2}}=(1-\alpha)^{3/2}\frac{\sqrt{|f'(\fraca\mu2)
|}}{f^2(\fraca\mu2)}>0,$
$s_{n,1}^2=\mathrm{O}(\sqrt{\beta_n}n^{-1})$ and $s_{n,2}^2=\mathrm{O}(\sqrt{\beta
_n}n^{-1})$. From this and (\ref{lalowerbound})
we obtain
%
%e16 ###
\begin{equation}\label{betauper}
\beta_n=\mathrm{O}(n^{-2/5}),
\end{equation}
with a similar bound being satisfied by $\gamma_n=\frac1 2-\tilde
{Q}_n(-\infty,0)$.

We turn now to the upper bound for $\mathcal{W}_2(\tilde{P}_n,P_0)$.
From the triangle inequality we get
\begin{eqnarray*}
\mathcal{W}_2(\tilde{P}_n,P_0)&\leq&\biggl (\int_0^{\fraca12}
(\tilde{F}_n^{-1}-F_0^{-1})^2  \biggr)^{1/2}+
 \biggl(\int_{\fraca12}^1 (\tilde{F}_n^{-1}-F_0^{-1})^2
\biggr)^{1/2} \\[-3pt]
&\leq&\mathcal{W}_2(\tilde{P}_n,\tilde{Q}_n)+ \biggl(\int
_0^{\fraca12} (\tilde{G}_n^{-1}-F_0^{-1})^2  \biggr)^{1/2}+
 \biggl(\int_{\fraca12}^1 (\tilde{F}_n^{-1}-F_0^{-1})^2  \biggr)^{1/2}.
\end{eqnarray*}
We consider next $\int_{\fraca1 2}^1 (\tilde{F}_n^{-1}-F_0^{-1})^2$.
Since $\tilde{P}_n\in\mathcal{R}_{\alpha_n}(P)$ we have
%
%e17 ###
\begin{equation}\label{ineq1}
\tilde{F}_n^{-1}(t)\leq-\frac\mu2 +F^{-1}\bigl(\alpha_n+(1-\alpha_n)t\bigr),
\qquad t\in(0,1).
\end{equation}
Keeping the above notation for $\beta_n$, assume first that $\beta
_n\leq0$.
Then
%
%e18 ###
\begin{equation}\label{ineq2}
\tilde{F}_n^{-1}(t)\geq-\frac\mu2 +F^{-1}\biggl(\alpha+(1-\alpha
_n)t+{\frac K {2\sqrt{n}}}\biggr), \qquad t\in\biggl({\frac1 2},1\biggr)
\end{equation}
(this follows upon noting that $\tilde{F}_n^{-1}({\frac1
2}+)\geq0$ and $\tilde{F}_n^{-1}(t)={F}^{-1}(h^{-1}(t))$, $h^{-1}$
growing with slope at least $1-\alpha_n$). For\vadjust{\goodbreak} $t
\in({\frac1 2},1)$, (\ref{ineq1}) and (\ref{ineq2}) still
hold if we replace $\tilde{F}_n^{-1}$ by~${F}_0^{-1}$. Hence,
in this case $\int_{\fraca1 2}^1 (\tilde
{F}_n^{-1}-F_0^{-1})^2\leq\int_{\fraca1 2}^1  (F^{-1}(\alpha
_n+(1-\alpha_n)t) -
F^{-1}(\alpha_n+(1-\alpha_n)t    -{\frac K
{2\sqrt{n}}})  )^2\,\mathrm{d}t=:s_{n,3}^2.$

If $\beta_n>0$, then, arguing as above, we have
%
%e19 ###
\begin{equation}\label{ineq3}
\tilde{F}_n^{-1}(t)\geq-\frac\mu2 +F^{-1}\biggl(\alpha+(1-\alpha_n)(t-\beta
_n)+{\frac K {2\sqrt{n}}}\biggr), \qquad t\in\biggl({\frac1
2}+\beta_n,1\biggr),
\end{equation}
while (\ref{ineq4}) holds in $(0,{\frac1 2}+\beta_n)$.
Now we use the bound
\[
\biggl(\int_{\fraca1 2}^1 (\tilde
{F}_n^{-1}-F_0^{-1})^2\biggr)^{1/2}\leq
\biggl(\int_{\fraca1 2}^{\fraca1 2+\beta_n} (\tilde{F_n}^{-1} -F_0^{-1})^2\biggr)^{1/2}+\biggl(\int_{\fraca1 2+\beta_n}^1 (\tilde
{F}_n^{-1}-F_0^{-1})^2\biggr)^{1/2}
\]
 and proceed as follows. For $t\in
({\frac1 2}+\beta_n,1)$
(\ref{ineq1}) and (\ref{ineq3}) hold again after replacing $\tilde
{F}_n^{-1}$ by $F_0^{-1}$. This and the triangle inequality yield
%
%e20 ###
\begin{eqnarray}\label{cotafinal1}\nonumber
 && \biggl(\int_{\fraca1 2+\beta_n}^1 (\tilde
{F}_n^{-1}-F_0^{-1})^2 \biggr)^{1/2} \hspace*{30pt}\\[-3pt]
&& \quad   \leq
 \biggl(\int_{\fraca1 2+\beta_n}^1\bigl(F^{-1}\bigl(\alpha+(1-\alpha
)t\bigr)-F^{-1}\bigl(\alpha+(1-\alpha)(t-\beta_n)\bigr)\bigr)^2\,\mathrm{d}t
\biggr)^{1/2}\nonumber\hspace*{30pt}
\\[-9pt]
\\[-9pt]
&& \qquad {} + 2 \biggl(\int_{\fraca1 2}^1\biggl(F^{-1}\bigl(\alpha_n+(1-\alpha
_n)t\bigr)-F^{-1}\biggl(\alpha_n+(1-\alpha_n)t-{\frac K {2\sqrt
{n}}}\biggr)\biggr)^2\,\mathrm{d}t \biggr)^{1/2}
\nonumber\hspace*{30pt}\\[-3pt]
 && \quad = \sqrt{g_2(\beta_n)}+2s_{n,3}.\hspace*{30pt}
\nonumber
\end{eqnarray}
For the interval $({\frac1 2},{\frac1 2}+\beta
_n)$ we write $\underline{G}^{-1}(t)={
\frac\mu2}+F^{-1}((1-\alpha_n)t)$ (the minimal quantile function in
$\mathcal{R}_{\alpha_n}(Q)$). Then
$(\int_{\fraca1 2}^{\fraca1 2+\beta_n} (\tilde
{F}_n^{-1}-F_0^{-1})^2)^{1/2}\leq
(\int_{\fraca1 2}^{\fraca1 2+\beta_n} (\tilde{F}_n^{-1}-\underline
{G}^{-1})^2)^{1/2}+(\int_{\fraca1 2}^{\fraca1 2+\beta_n} (\underline
{G}^{-1}-F_0^{-1})^2)^{1/2}$.
We observe now that $\tilde{G}^{-1}(t)\geq\underline{G}_n^{-1}(t)$ and
also that, for $t \in({\frac1 2},{\frac1
2}+\beta_n)$, ${-\frac\mu2}+F^{-1}(\alpha+(1-\alpha
)(t-\beta_n))\leq0 \leq
{\frac\mu2}+F^{-1}((1-\alpha)t)$. Combining these facts
with (\ref{ineq4}) we obtain
\begin{eqnarray*}
|\tilde{F}_n^{-1}(t)-\underline{G}^{-1}(t)|&\leq&|\tilde
{F}_n^{-1}(t)-\tilde{G}_n^{-1}(t)|  \\
&&{}+\bigl| F^{-1}\bigl((1-\alpha_n)t\bigr)- F^{-1}\bigl((1-\alpha)t\bigr)\bigr|\\
&&{}+\biggl| F^{-1}\biggl(\alpha+(1-\alpha_n)(t-\beta_n)+{\frac K {2\sqrt
{n}}}\biggr)- F^{-1}\bigl(\alpha+(1-\alpha)(t-\beta_n)\bigr)\biggr|.
\end{eqnarray*}
As a consequence,
%
%e21 ###
\begin{eqnarray*}
 &&\biggl(\int_{\fraca1 2}^{\fraca1 2+\beta_n} (\tilde
{F}_n^{-1}-F_0^{-1})^2 \biggr)^{1/2}\\[-3pt]
&& \quad \leq\mathcal{W}_2 (\tilde{P}_n,\tilde
{Q}_n)+
 \biggl(\int_{\fraca1 2}^{\fraca1 2+\beta_n} \bigl(\mu+{F}^{-1}\bigl((1-\alpha
)t\bigr)-F^{-1}\bigl(\alpha+(1-\alpha)t\bigr)\bigr)^2\,\mathrm{d}t \biggr)^{1/2}
\\[-3pt]
&& \qquad {}+
2 \biggl(\int_{\fraca1 2}^{\fraca1 2+\beta_n}\bigl(F^{-1}\bigl((1-\alpha
_n)t\bigr)-F^{-1}\bigl((1-\alpha)t\bigr)\bigr)^2\,\mathrm{d}t \biggr)^{1/2}\\[-3pt]
&& \qquad {}+
 \biggl(\int_{\fraca1 2}^{\fraca1 2+\beta_n}\biggl(F^{-1}\biggl(\alpha+(1-\alpha
_n)(t-\beta_n)+{\frac K {2\sqrt{n}}}\biggr)\\[-3pt]
&&\phantom{\qquad {}+
 \biggl(\int_{\fraca1 2}^{\fraca1 2+\beta_n}\biggl(}{}-F^{-1}\bigl(\alpha+(1-\alpha
)(t-\beta_n)\bigr)\biggr)^2\,\mathrm{d}t \biggr)^{1/2}\\[-3pt]
&& \quad =
\mathcal{W}_2 (\tilde{P}_n,\tilde{Q}_n)+\sqrt{g_3(\beta_n)}+2
s_{n,4}+s_{n,5},\vspace*{-2pt}
\end{eqnarray*}
where $g_3(\beta)=\int_{\fraca1 2}^{\fraca1 2+\beta} (\mu
+{F}^{-1}((1-\alpha)t)-F^{-1}(\alpha+(1-\alpha)t))^2\,\mathrm{d}t$. Again a Taylor
expansion shows that
$g_3(\beta_n)=\mathrm{O}(\beta_n^3)=\mathrm{o}(n^{-1})$. Similarly, we get
$s_{n,j}=\mathrm{o}(n^{-1})$, $j=4,5$, and, as a consequence\vspace*{-2pt}
%
%e22 ###
\begin{equation}\label{trozodeenmedio}
 \biggl(\int_{\fraca1 2}^{\fraca1 2+\beta_n} (\tilde
{F}_n^{-1}-F_0^{-1})^2 \biggr)^{1/2}=\mathrm{O}(n^{-1/2}).\vspace*{-2pt}
\end{equation}
Collecting the estimates in (\ref{cotafinal1}) and (\ref{trozodeenmedio}),
we obtain\vspace*{-2pt}
%
%e23 ###
\begin{equation}\label{yaveremos}
 \biggl(\int_{\fraca1 2}^{1} (\tilde{F}_n^{-1}-F_0^{-1})^2
\biggr)^{1/2}\leq\sqrt{g_2(\beta_n)}+2s_{n,3}+\mathrm{O}(n^{-1/2}).\vspace*{-2pt}
\end{equation}
We note next that $F^{-1}$ has a bounded derivative and, as a consequence,
$s_{n,3}^2=\mathrm{O}(n^{-1})$. Similarly, we find that $g_2(\beta_n)=\mathrm{O}(\beta
_n^2)$. Summarizing,\vspace*{-2pt}
\[
 \biggl(\int_{\fraca1 2}^{1} (\tilde{F}_n^{-1}-F_0^{-1})^2
\biggr)^{1/2}=\mathrm{O}(n^{-\fraca2 5}).\vspace*{-2pt}
\]
A similar analysis works for $\int_0^{\fraca1 2} (\tilde
{G}_n^{-1}-F_0^{-1})^2$ and completes the proof.\vspace*{-2pt}
\end{pf*}

\begin{pf*}{Proof of Proposition \ref{empiricalcost}}
We take $(X_{1,1},Y_{1,1})$ to be an optimal coupling for $P$ and $Q$
with respect to the $\|x-y\|^p$-cost and $(X_{1,i},Y_{1,i})$, $2\leq
i\leq n$, and $(X_{2,j},Y_{2,j})$, $1\leq j\leq m$, independent copies
of $(X_{1,1},Y_{1,1})$ (hence $E\|X_{i,j} -Y_{i,j}\|^p=\mathcal{W}_p^p
(P,Q)$). Then
$S_{n,m}=\min_\pi(a(\pi))^{1/p}$ and $T_{n,m}=\min_\pi(b(\pi
))^{1/p}$, where\vspace*{-2pt}
\[
a(\pi)=\sum_{1\leq i\leq n,1\leq j \leq m}\pi_{i,j}\|X_{1,i}-X_{2,j} \|^p,\vspace*{-2pt}
\]
$b(\pi)$ is defined similarly by replacing $X_{i,j}$ by $Y_{i,j}$ and
$\pi$ takes values in the set of $n\times m$ matrices with non-negative
entries $\pi_{i,j}$ such that
$\sum_{1\leq j\leq m} \pi_{i,j}=\frac1 n$ and $\sum_{1\leq i\leq n} \pi
_{i,j}=\frac1 m$.

We observe next that, by the triangle inequality,\vspace*{-2pt}
\begin{eqnarray*}
|a(\pi)^{1/p}-b(\pi)^{1/p}|&\leq&  \biggl(\sum_{1\leq i\leq
n,1\leq j \leq m}\pi_{i,j}\|(X_{1,i}-X_{2,j})-(Y_{1,i}-Y_{2,j}) \|^p
 \biggr)^{1/p}
 \\[-3pt]
&\leq&  \biggl(\frac1 n \sum_{1\leq i\leq n}\|X_{1,i}-Y_{1,i}\|^p  \biggr)^{1/p}+
 \biggl(\frac1 m \sum_{1\leq j\leq m}\|X_{2,j}-Y_{2,j}\|^p
 \biggr)^{1/p}.\vspace*{-2pt}\vadjust{\goodbreak}
\end{eqnarray*}
As a consequence, we have that $|S_{n,m}-T_{n,m}| $ is upper bounded by
the right-hand side of the above display and, from the elementary inequality
$(a+b)^p\leq2^{p-1}a^p +2^{p-1}b^p$ for non-negative $a,b$, we get
\begin{eqnarray*}
E(S_{n,m}-T_{n,m})^p&\leq& 2^{p-1}E\|X_{1,1} -Y_{1,1} \|^p + 2^{p-1}E\|
X_{2,1} -Y_{2,1} \|^p\\
&=& 2^p \mathcal{W}_p^p (P,Q).
\end{eqnarray*}
This completes the proof.
\end{pf*}
\end{appendix}

\section*{Acknowledgements}
Research partially supported by the Spanish Ministerio de Ciencia e
Innovaci\'on,  Grant MTM2008-06067-C02-01,
and 02 and by the Consejer\'{\i}a de Educaci\'on y Cultura de la Junta
de Castilla y Le\'on, GR150.

The authors would like to thank two anonymous referees for their
careful reading of the manuscript,
their suggestions and the pointers to relevant references that helped
us to greatly
improve our original version.

% imsref loaded by smiklovaite, 2011-10-05 07:08:54
% imsref loaded by smiklovaite, 2011-10-05 07:14:04

\printhistory


\begin{thebibliography}{24}
% BibTex style file: bej.bst, 2010-01-21
% Default style options (sort=1,type=number).
% Used options (sort=1,type=number).

%b1 ###
\bibitem{Pedro}
\begin{barticle}[mr]
\bauthor{\bsnm{{\'A}lvarez-Esteban},~\bfnm{Pedro~C{\'e}sar}\binits{P.C.}},
  \bauthor{\bparticle{del} \bsnm{Barrio},~\bfnm{Eustasio}\binits{E.}},
  \bauthor{\bsnm{Cuesta-Albertos},~\bfnm{Juan~Antonio}\binits{J.A.}} \AND
  \bauthor{\bsnm{Matr{\'a}n},~\bfnm{Carlos}\binits{C.}}
(\byear{2008}).
\btitle{Trimmed comparison of distributions}.
\bjournal{J. Amer. Statist. Assoc.}
\bvolume{103}
\bpages{697--704}.
\bid{doi={10.1198/016214508000000274}, issn={0162-1459}, mr={2435470}}
\end{barticle}
\endbibitem

%b2 ###
\bibitem{Pedro2}
\begin{barticle}[auto:STB|2011/09/12|07:03:23]
\bauthor{\bsnm{{\'A}lvarez-Esteban},~\bfnm{P.~C.}\binits{P.C.}},
  \bauthor{\bparticle{del} \bsnm{Barrio},~\bfnm{E.}\binits{E.}},
  \bauthor{\bsnm{Cuesta-Albertos},~\bfnm{J.~A.}\binits{J.A.}} \AND
  \bauthor{\bsnm{Matr{\'a}n},~\bfnm{C.}\binits{C.}}
(\byear{2011}).
\btitle{Uniqueness and approximate computation of optimal incomplete
  transportation plans}. \bjournal{Ann. Inst. H. Poincar\'e Probab.
  Statist.} \bvolume{47} \bpages{358--375}.
  \bid{mr={2814414}}
\end{barticle}
\endbibitem

%b3 ###
\bibitem{Pedro3}
\begin{barticle}[auto:STB|2011/09/12|07:03:23]
\bauthor{\bsnm{{\'A}lvarez-Esteban},~\bfnm{P.~C.}\binits{P.C.}},
  \bauthor{\bparticle{del} \bsnm{Barrio},~\bfnm{E.}\binits{E.}},
  \bauthor{\bsnm{Cuesta-Albertos},~\bfnm{J.~A.}\binits{J.A.}} \AND
  \bauthor{\bsnm{Matr{\'a}n},~\bfnm{C.}\binits{C.}}
(\byear{2010}).
\btitle{Assessing when a sample is mostly normal}.
\bjournal{Comput. Statist. Data Anal.}
\bvolume{54}
\bpages{2914--2925}.
\end{barticle}
\endbibitem

%b4 ###
\bibitem{Bickel}
\begin{barticle}[mr]
\bauthor{\bsnm{Bickel},~\bfnm{Peter~J.}\binits{P.J.}} \AND
  \bauthor{\bsnm{Freedman},~\bfnm{David~A.}\binits{D.A.}}
(\byear{1981}).
\btitle{Some asymptotic theory for the bootstrap}.
\bjournal{Ann. Statist.}
\bvolume{9}
\bpages{1196--1217}.
\bid{issn={0090-5364}, mr={0630103}}
\end{barticle}
\endbibitem

%b5 ###
\bibitem{Buja}
\begin{barticle}[mr]
\bauthor{\bsnm{Buja},~\bfnm{Andreas}\binits{A.}}
(\byear{1986}).
\btitle{On the {H}uber--{S}trassen theorem}.
\bjournal{Probab. Theory Relat. Fields}
\bvolume{73}
\bpages{149--152}.
\bid{doi={10.1007/BF01845998}, issn={0178-8051}, mr={0849070}}
\end{barticle}
\endbibitem

%b6 ###
\bibitem{CaffarelliMcCann}
\begin{barticle}[mr]
\bauthor{\bsnm{Caffarelli},~\bfnm{Luis~A.}\binits{L.A.}} \AND
  \bauthor{\bsnm{McCann},~\bfnm{Robert~J.}\binits{R.J.}}
(\byear{2010}).
\btitle{Free boundaries in optimal transport and {M}onge--{A}mp\`ere obstacle
  problems}.
\bjournal{Ann. of Math. (2)}
\bvolume{171}
\bpages{673--730}.
\bid{doi={10.4007/annals.2010.171.673}, issn={0003-486X}, mr={2630054}}
\end{barticle}
\endbibitem

%b7 ###
\bibitem{Cascos2}
\begin{barticle}[mr]
\bauthor{\bsnm{Cascos},~\bfnm{Ignacio}\binits{I.}} \AND
  \bauthor{\bsnm{L{\'o}pez-D{\'{\i}}az},~\bfnm{Miguel}\binits{M.}}
(\byear{2008}).
\btitle{Consistency of the {$\alpha$}-trimming of a probability. {A}pplications
  to central regions}.
\bjournal{Bernoulli}
\bvolume{14}
\bpages{580--592}.
\bid{doi={10.3150/07-BEJ109}, issn={1350-7265}, mr={2544103}}
\end{barticle}
\endbibitem

%b8 ###
\bibitem{CsH93}
\begin{bbook}[mr]
\bauthor{\bsnm{Cs{\"o}rg{\H{o}}},~\bfnm{Mikl{\'o}s}\binits{M.}} \AND
  \bauthor{\bsnm{Horv{\'a}th},~\bfnm{Lajos}\binits{L.}}
(\byear{1993}).
\btitle{Weighted Approximations in Probability and Statistics}.
\bseries{Wiley Series in Probability and Mathematical Statistics: Probability
  and Mathematical Statistics}.
\baddress{Chichester}: \bpublisher{Wiley}.
\bid{mr={1215046}}
\end{bbook}
\endbibitem

%b9 ###
\bibitem{Cuesta2}
\begin{barticle}[mr]
\bauthor{\bsnm{Cuesta-Albertos},~\bfnm{J.~A.}\binits{J.A.}},
  \bauthor{\bsnm{Gordaliza},~\bfnm{A.}\binits{A.}} \AND
  \bauthor{\bsnm{Matr{\'a}n},~\bfnm{C.}\binits{C.}}
(\byear{1997}).
\btitle{Trimmed {$k$}-means: An attempt to robustify quantizers}.
\bjournal{Ann. Statist.}
\bvolume{25}
\bpages{553--576}.
\bid{doi={10.1214/aos/1031833664}, issn={0090-5364}, mr={1439314}}
\end{barticle}
\endbibitem

%b10 ###
\bibitem{Barrio}
\begin{barticle}[mr]
\bauthor{\bparticle{del} \bsnm{Barrio},~\bfnm{Eustasio}\binits{E.}},
  \bauthor{\bsnm{Gin{\'e}},~\bfnm{Evarist}\binits{E.}} \AND
  \bauthor{\bsnm{Utzet},~\bfnm{Frederic}\binits{F.}}
(\byear{2005}).
\btitle{Asymptotics for {$L\sb 2$} functionals of the empirical quantile
  process, with applications to tests of fit based on weighted {W}asserstein
  distances}.
\bjournal{Bernoulli}
\bvolume{11}
\bpages{131--189}.
\bid{doi={10.3150/bj/1110228245}, issn={1350-7265}, mr={2121458}}
\end{barticle}
\endbibitem

%b11 ###
\bibitem{Figalli}
\begin{barticle}[mr]
\bauthor{\bsnm{Figalli},~\bfnm{Alessio}\binits{A.}}
(\byear{2010}).
\btitle{The optimal partial transport problem}.
\bjournal{Arch. Ration. Mech. Anal.}
\bvolume{195}
\bpages{533--560}.
\bid{doi={10.1007/s00205-008-0212-7}, issn={0003-9527}, mr={2592287}}
\end{barticle}
\endbibitem

%b12 ###
\bibitem{Garcia}
\begin{barticle}[mr]
\bauthor{\bsnm{Garc{\'{\i}}a-Escudero},~\bfnm{Luis~A.}\binits{L.A.}},
  \bauthor{\bsnm{Gordaliza},~\bfnm{Alfonso}\binits{A.}},
  \bauthor{\bsnm{Matr{\'a}n},~\bfnm{Carlos}\binits{C.}} \AND
  \bauthor{\bsnm{Mayo-Iscar},~\bfnm{Agustin}\binits{A.}}
(\byear{2008}).
\btitle{A general trimming approach to robust cluster analysis}.
\bjournal{Ann. Statist.}
\bvolume{36}
\bpages{1324--1345}.
\bid{doi={10.1214/07-AOS515}, issn={0090-5364}, mr={2418659}}
\end{barticle}
\endbibitem

%b13 ###
\bibitem{Gorda}
\begin{barticle}[mr]
\bauthor{\bsnm{Gordaliza},~\bfnm{Alfonso}\binits{A.}}
(\byear{1991}).
\btitle{Best approximations to random variables based on trimming procedures}.
\bjournal{J.~Approx. Theory}
\bvolume{64}
\bpages{162--180}.
\bid{doi={10.1016/0021-9045(91)90072-I}, issn={0021-9045}, mr={1091467}}
\end{barticle}
\endbibitem

%b14 ###
\bibitem{Gower}
\begin{bincollection}[auto:STB|2011/09/12|07:03:23]
\bauthor{\bsnm{Gower},~\bfnm{J.~C.}\binits{J.C.}}
(\byear{2006}).
\btitle{Measures of similarity, dissimilarity, and distance}.
In \bbooktitle{Encyclopedia of Statistical Sciences}, \bedition{2nd ed.}
(\beditor{\bfnm{Samuel}\binits{S.}~\bsnm{Kotz}},
  \beditor{\bfnm{Campbell~B.}\binits{C.B.}~\bsnm{Read}},
  \beditor{\bfnm{N.}\binits{N.}~\bsnm{Balakrishnan}} \AND
  \beditor{\bfnm{Brani}\binits{B.}~\bsnm{Vidakovic}}, eds.)
\bvolume{12}
\bpages{7730--7738}.
\baddress{New York}: \bpublisher{Wiley}.
\end{bincollection}
\endbibitem

%b15 ###
\bibitem{Huber64}
\begin{barticle}[mr]
\bauthor{\bsnm{Huber},~\bfnm{Peter~J.}\binits{P.J.}}
(\byear{1964}).
\btitle{Robust estimation of a location parameter}.
\bjournal{Ann. Math. Statist.}
\bvolume{35}
\bpages{73--101}.
\bid{issn={0003-4851}, mr={0161415}}
\end{barticle}
\endbibitem

%b16 ###
\bibitem{Huber65}
\begin{barticle}[mr]
\bauthor{\bsnm{Huber},~\bfnm{Peter~J.}\binits{P.J.}}
(\byear{1965}).
\btitle{A robust version of the probability ratio test}.
\bjournal{Ann. Math. Statist.}
\bvolume{36}
\bpages{1753--1758}.
\bid{issn={0003-4851}, mr={0185747}}
\end{barticle}
\endbibitem

%b17 ###
\bibitem{HuberStrassen}
\begin{barticle}[mr]
\bauthor{\bsnm{Huber},~\bfnm{Peter~J.}\binits{P.J.}} \AND
  \bauthor{\bsnm{Strassen},~\bfnm{Volker}\binits{V.}}
(\byear{1973}).
\btitle{Minimax tests and the {N}eyman--{P}earson lemma for capacities}.
\bjournal{Ann. Statist.}
\bvolume{1}
\bpages{251--263}.
\bid{issn={0090-5364}, mr={0356306}}
\end{barticle}
\endbibitem

%b18 ###
\bibitem{Maronna2}
\begin{barticle}[mr]
\bauthor{\bsnm{Maronna},~\bfnm{Ricardo}\binits{R.}}
(\byear{2005}).
\btitle{Principal components and orthogonal regression based on robust scales}.
\bjournal{Technometrics}
\bvolume{47}
\bpages{264--273}.
\bid{doi={10.1198/004017005000000166}, issn={0040-1706}, mr={2164700}}
\end{barticle}
\endbibitem

%b19 ###
\bibitem{Camblor}
\begin{barticle}[mr]
\bauthor{\bsnm{Mart{\'{\i}}nez-Camblor},~\bfnm{P.}\binits{P.}},
  \bauthor{\bsnm{De~U{\~n}a-{\'A}lvarez},~\bfnm{J.}\binits{J.}} \AND
  \bauthor{\bsnm{Corral},~\bfnm{N.}\binits{N.}}
(\byear{2008}).
\btitle{{$k$}-sample test based on the common area of kernel density
  estimators}.
\bjournal{J. Statist. Plann. Inference}
\bvolume{138}
\bpages{4006--4020}.
\bid{doi={10.1016/j.jspi.2008.02.008}, issn={0378-3758}, mr={2455983}}
\end{barticle}
\endbibitem

%b20 ###
\bibitem{Massart}
\begin{bbook}[mr]
\bauthor{\bsnm{Massart},~\bfnm{Pascal}\binits{P.}}
(\byear{2007}).
\btitle{Concentration Inequalities and Model Selection}.
\bseries{Lecture Notes in Math.}
\bvolume{1896}.
\baddress{Berlin}: \bpublisher{Springer}.
%  Saint-Flour, July 6--23, 2003, With a foreword by Jean Picard}.
\bid{mr={2319879}}
\end{bbook}
\endbibitem

%b21 ###
\bibitem{Raghavachari1973}
\begin{barticle}[mr]
\bauthor{\bsnm{Raghavachari},~\bfnm{M.}\binits{M.}}
(\byear{1973}).
\btitle{Limiting distributions of {K}olmogorov--{S}mirnov type statistics under
  the alternative}.
\bjournal{Ann. Statist.}
\bvolume{1}
\bpages{67--73}.
\bid{issn={0090-5364}, mr={0346976}}
\end{barticle}
\endbibitem

%b22 ###
\bibitem{Rieder}
\begin{barticle}[mr]
\bauthor{\bsnm{Rieder},~\bfnm{Helmut}\binits{H.}}
(\byear{1977}).
\btitle{Least favorable pairs for special capacities}.
\bjournal{Ann. Statist.}
\bvolume{5}
\bpages{909--921}.
\bid{issn={0090-5364}, mr={0468005}}
\end{barticle}
\endbibitem

%b23 ###
\bibitem{Rousseeuw1}
\begin{bincollection}[mr]
\bauthor{\bsnm{Rousseeuw},~\bfnm{Peter}\binits{P.}}
(\byear{1985}).
\btitle{Multivariate estimation with high breakdown point}.
In \bbooktitle{Mathematical Statistics and Applications, {V}ol. {B} ({B}ad
  {T}atzmannsdorf, 1983)}
  (\beditor{W.~Grossmann}, \beditor{G.~Pflug}, \beditor{I.~Vincze} and \beditor{W.~Werz}, eds.)
\bpages{283--297}.
\baddress{Dordrecht}: \bpublisher{Reidel}.
\bid{mr={0851060}}
\end{bincollection}
\endbibitem

\end{thebibliography}
\end{document}